\documentstyle[leqno]{article}
\begin{document}
\renewcommand{\theequation}{\arabic{section}.\arabic{equation}}
\newcommand{\be}{\begin{eqnarray}}
\newcommand{\en}{\end{eqnarray}}
\newcommand{\no}{\nonumber}
\newcommand{\la}{\lambda}
\newcommand{\laa}{\Lambda}
\newcommand{\ep}{\epsilon}
\newcommand{\de}{\delta}
\newcommand{\pl}{\parallel}
\newcommand{\ov}{\overline}
\newcommand{\bet}{\beta}
\newcommand{\al}{\alpha}
\newcommand{\fr}{\frac}
\newcommand{\pa}{\partial}
\newcommand{\lan}{\langle}
\newcommand{\ra}{\rangle}
\newcommand{\we}{\wedge}
\newcommand{\om}{\Omega}
\newcommand{\na}{\nabla}
\newcommand{\D}{\Delta}
\newcommand{\vs}{\vskip0.3cm}
\renewcommand{\thefootnote}{}
\title {Universal Bounds for Eigenvalues of the Polyharmonic Operators  }
\footnotetext{
\hspace*{2ex}Key words and phrases: Universal bounds,  eigenvalues, polyharmonic operator, Riemannian manifods, Euclidean space, spheres.}
\author{J\"urgen Jost, Xianqing Li-Jost, Qiaoling Wang, Changyu Xia }\date{} \maketitle
\begin{abstract}
We study   eigenvalues of  polyharmonic
operators on compact Riemannian manifolds with boundary (possibly
empty). In particular, we prove a universal inequality for  the  eigenvalues of
the polyharmonic operators on compact domains in a Euclidean space. This inequality controls the $k$th eigenvalue by the lower eigenvalues, independently of the particular geometry of the domain.
Our inequality is  sharper than the known Payne-P\'olya-Weinberg type
inequality and  also covers the important Yang inequality on
eigenvalues of the Dirichlet Laplacian. We also prove universal  inequalities
for the lower order eigenvalues of the polyharmonic operator on compact
domains in a Euclidean space which in the case of the biharmonic
operator and the buckling problem strengthen the estimates obtained by Ashbaugh.
Finally,   we  prove universal
inequalities  for eigenvalues of  polyharmonic
operators of any order on compact domains in the sphere.
\end{abstract}


\section{Introduction}
Let $\om$ be a connected bounded domain with smooth boundary in  an $n(\geq 2)$-dimensional Euclidean space  ${I\!\!R}^n$ and let $\nu$ be the outward unit normal vector field of $\pa \om$.
Denote by $\Delta $ the Laplace operator on ${I\!\!R}^n$ and let  $l$ be a positive integer. Solutions of $\Delta u=0$ on a domain $\om \subset {I\!\!R}^n$ are of course the classical harmonic functions which describe the equilibrium position of an elastic homogeneous membrane. Solutions of $\Delta^2 u=0$ are called biharmonic, and they model equilibria of homogeneous plates. Similarly, solutions of  $\Delta^l u=0,\ l\in {I\!\!N},$ are called polyharmonic.

One then naturally considers the
eigenvalue problem
\be
& & (-\Delta)^l u=\la u\ \ \ {\rm in} \ \ \om, \\
\no  & &  u|_{\pa \om}=\left. \fr{\pa u}{\pa \nu}\right|_{\pa \om}=\cdots =\left. \fr{\pa^{l-1} u}{\pa \nu^{l-1}}\right|_{\pa \om}=0.
\en
Let
\be\no
 0<\la_1\leq \la_2\leq\la_3\leq\cdots,
\en denote the successive eigenvalues, where
each eigenvalue is repeated according to its multiplicity.\\
The case $l=1$ is of course well-studied, since the work of Weyl [We] and Courant-Hilbert [CH].  But also for $l\ge 2$, polyharmonic functions have interesting applications in physics. For
example, the Airy function in mechanics is a bi-harmonic function. More generally, a clamped plate in equilibrium is a solution of the bi-harmonic problem
\be\no
& & (-\Delta)^2 v=0\ \ \ {\rm in} \ \ \om, \\
\no  & &  v|_{\pa \om}=\left. \fr{\pa v}{\pa \nu}\right|_{\pa \om}=0
\en
in a two-dimensional $\om$. An oscillating clamped plate then satisfies
\be\no
& & (-\Delta)^2 v +v_{tt}=0\ \ \ {\rm in} \ \ \om, \ t\ge 0,\\
\no  & &  v|_{\pa \om}=\left. \fr{\pa v}{\pa \nu}\right|_{\pa \om}=0,
\en
and a separation of variables $v(x,y,t)=u(x,y)g(t)$ leads to the eigenvalue problem
\be\no
& & (-\Delta)^2 u=\la u\ \ \ {\rm in} \ \ \om, \\
\no  & &  u|_{\pa \om}=\left. \fr{\pa u}{\pa \nu}\right|_{\pa \om}=0.
\en
This problem has been studied already by Courant [Co]. He derived the Weyl type law
\be\no
\la_k \sim \left( \frac{4\pi k}{area(\om)}\right)^2.
\en

In this paper, we investigate the  eigenvalues of the problem
(1.1) for general $l$. We are interested in so-called universal properties, that is, properties that do not depend on the specific domain $\om$, but only on its dimension $n$. These universal properties then take the form of relations between different eigenvalues. Naturally, the first eigenvalue $\la_1$ plays a distinguished role. Since this eigenvalue can often be estimated in terms of the geometry of $\om$, one can then also derive geometric estimates for higher eigenvalues from such universal bounds, but this is not explored in the present paper.

Let us now put our results into the context of those known for $l=1$.  Payne, P\'olya and Weinberger     proved in [PPW1] and [PPW2] that
$$\fr{\la_2}{\la_1}\leq 3 \ \ \ {\rm for\ \ \ }  \om\subset {I\!\!R}^2$$
and conjectured that
$$\left. \fr{\la_2}{\la_1}\leq \fr{\la_2}{\la_1}\right|_{\rm disk}$$
with equality if and only if $\om$ is a disk. For $n\geq 2$, the analogous statements are
$$\fr{\la_2}{\la_1}\leq 1+\fr 4n  \ \ \ {\rm for\ \ \ }  \om\subset {I\!\!R}^n,$$
and the {\it PPW conjecture}
$$\left. \fr{\la_2}{\la_1}\leq \fr{\la_2}{\la_1}\right|_{n-{\rm ball}},$$
with equality if and only if $\om$ is an $n$-ball.
This important {\it PPW conjecture} was solved by Ashbaugh and Benguria in  [AB1], [AB2], [AB3]. In [PPW2], Payne, P\'olya and Weinberger also proved the bound
\be
\la_{k+1}-\la_k\leq \fr 2k\sum_{i=1}^k\la_i, \ \ k=1, 2,\cdots,
\en
for $\om\subset {I\!\!R}^2$. This result easily extends to $\om\subset {I\!\!R}^n$ as
\be
\la_{k+1}-\la_k\leq \fr 4{kn}\sum_{i=1}^k\la_i, \ \ k=1, 2,\cdots,
\en
Many interesting works have been done in generalizing (1.3), e. g., in [A1], [A2], [AH], [CY1], [Ha], [HM1], [HM2], [HP], [HS], [HY], [LeP], [Y]. Here we  mention two results in this direction.
In 1980, Hile and Protter proved [HP]
\be
\sum_{i=1}^k\fr{\la_i}{\la_{k+1}-\la_i}\geq \fr{kn}4, \ \ \ {\rm for} \ k=1,2,\cdots.
\en
In 1991, Yang [Y] proved the following much stronger inequality:
\be
\sum_{i=1}^k(\la_{k+1}-\la_i)\left(\la_{k+1}-\left(1+\fr 4n\right)\la_i\right)\leq 0, \ \ \
{\rm for \ } k=1, 2,\cdots .
\en
The inequality (1.5), as observed by Yang himself, and as later proved, e. g., in [A1], [A2], [AH],  is the strongest of the classical inequalities that are derived following the scheme devised by Payne-P\'olya-Weinberger. Yang's inequality provided a marked improvement for eigenvalues of large index. Recently, some Yang type inequalities on eigenvalues of the problem (1.1) for the case $l>1$ have been proved in [CY2], [CY3], [WX1], [WX2] and [WC]. We remark that there is an error in the line below (3.1) of [WX1] where boundary terms are dropped from an integration by parts, without a reason that these terms should vanish.

For  general $l$, the Payne-P\'olya-Weinberg type inequality reads (cf. [CQ], [H]):
\be
\la_{k+1}\leq \la_k+\fr{4l(n+2l-2)}{n^2k^2}\left(\sum_{i=1}^k\la_i^{1/l}\right)\left(\sum_{i=1}^k\la_i^{(l-1)/l}\right),
\en
In this paper, we obtain a universal inequality  of Yang type for the  eigenvalues  of the problem (1.1)  for any $l$. Indeed, we  consider the more general eigenvalue problem:
\be
& & (-\Delta)^l u=\la u\ \ \ {\rm in} \ \ M,
\\ \no
& & u|_{\pa M}=\left. \fr{\pa u}{\pa \nu}\right|_{\pa M}=\cdots =\left. \fr{\pa^{l-1} u}{\pa \nu^{l-1}}\right|_{\pa M}=0,
\en
where $M$ is a compact Riemannian manifold with boundary (possibly empty), $\Delta$ is the Laplacian operator on
$M$ (for general results for the case $l=1$, see e.g. [Ch]). We will prove  a general inequality for the eigenvalues of the problem (1.7) (see Theorem 2.1). By using this inequality, we  show that when $M$ is a bounded connected domain in ${I\!\!R}^n$ with smooth boundary then the eigenvalues of the problem (1.1) satisfy (see Theorem 3.1):
\be
& &
\sum_{i=1}^k\left(\la_{k+1}-\la_i\right)^2
\\ \no &\leq& \left(\frac{4l(n+2l-2)}{n^2}\right)^{1/2}
\left(\sum_{i=1}^k\left(\la_{k+1}-\la_i\right)^2\la_i^{(l-1)/l}\right)^{1/2}
\left(\sum_{i=1}^k\left(\la_{k+1}-\la_i\right)\la_i^{1/l}\right)^{1/2}.
\en When $l=1$, (1.8) is just Yang's inequality (1.5). As a
consequence of (1.8), we have the following two estimates for the
$(k+1)$-th eigenvalue in terms of the first $k$-eigenvalues of the
problem (1.1):
\be & & \la_{k+1}\\ \no &\leq &\left\{\left(\fr{2l(n+2l-2)}{k^2n^2}\right)^2\left(\sum_{i=1}^k\la_i^{(l-1)/l}
\right)^2\left(\sum_{i=1}^k\la_i^{1/l}\right)^2-\fr
1k\sum_{i=1}^k\left(\la_i-\fr 1k\sum_{j=1}^k \la_j
\right)^2\right\}^{1/2} \\ \no & &
+\fr
1k\sum_{i=1}^k\la_i+\fr{2l(n+2l-2)}{k^2n^2}\left(\sum_{i=1}^k\la_i^{(l-1)/l}
\right)\left(\sum_{i=1}^k\la_i^{1/l}\right)
\en

\be & & \la_{k+1}\\ \no &\leq
&\left\{\left(\fr{2l(n+2l-2)}{n^2}\fr
1k\sum_{i=1}^k\la_i\right)^2-\left(1+\fr{4l(n+2l-2)}{n^2}\right)\fr
1k\sum_{j=1}^k\left(\la_j-\fr
1k\sum_{i=1}^k\la_i\right)^2\right\}^{1/2}\\ \no & &
+\left(1+\fr{2l(n+2l-2)}{n^2}\right)\fr 1k\sum_{i=1}^k\la_i. \en
Notice that (1.9) is much stronger than (1.6).

In [AB4], Ashbaugh and Benguria showed that when $l=1$, the first $n+1$ eigenvalues of
the problem (1.1) satisfy the inequality
\be \la_2+\la_3+\cdots
+\la_{n+1}\leq (n+4)\la_1. \en Ashbaugh showed in [A1] that when
$l=2$,
 \be \la_2+\la_3+\cdots +\la_{n+1}\leq (n+24)\la_1.
\en
In this paper, we prove a similar inequality for any $l$ which
covers the inequality (1.11) when $l=1$ and improves (1.12) when
$l=2$ (Cf. Theorem 4.1). The reason why the dimension $n$ comes in here is that  coordinate functions in $n$-dimensional Euclidean space when used as test functions yield useful inequalities.

Consider now the so called buckling problem :
\be \Delta^2 u=-\laa \Delta u \ \ \ {\rm in} \ \ \om ,  \ \
\ \ \left. u\right|_{\pa \om}=\left. \fr{\pa u}{\pa \nu}\right|_{\pa
\om}=0, \en
where $\om $ is a bounded connected domain in $I\!\!R^n$.

Let \be\no
 0<\laa_1\leq\laa_2\leq \laa_{3}\leq \cdots
\en denote the successive eigenvalues for (1.13).
Payne, P\'olya and Weinberger [PPW2] proved
$$\laa_2/\laa_1<3 \ \ \ \ {\rm for }\ \om\subset I\!\!R^2.$$
For $\om\subset I\!\!R^n$ this reads
$$\laa_2/\laa_1<1+4/n.$$
Subsequently, Hile and Yeh [HY]  obtained the improved bound
$$
\fr{\laa_2}{\laa_1}\leq \fr{n^2+8n+20}{(n+2)^2}\ \ \ \ \ {\rm for \ } \om\subset I\!\!R^n.
$$
Ashbaugh [A] proved :
\be \sum_{i=1}^n\laa_{i+1}\leq(n+4)\laa_1. \en
 Cheng and Yang  [CY2] obtained: \be \sum_{i=1}^k(\laa_{k+1}-\laa_i)^2\leq
\fr{4(n+2)}{n^2}\sum_{i=1}^k(\laa_{k+1}-\laa_i)\laa_i. \en
In this paper, we will prove the following inequality which strengthens (1.14) (Cf. Theorem 4.2):
\be \sum_{i=1}^n\laa_{i+1}
+\fr{4(\laa_2-\laa_1)}{n+4}\leq (n+4)\laa_1.\en
We will also show that the first $(n+1)$ eigenvalues of the following more general problem
\be
& & (-\Delta)^l u=-\laa\Delta u\ \ \ {\rm in} \ \ \om, \\
\no  & &  u|_{\pa \om}=\left. \fr{\pa u}{\pa \nu}\right|_{\pa
\om}=\cdots =\left. \fr{\pa^{l-1} u}{\pa \nu^{l-1}}\right|_{\pa
\om}=0, \en where $l\geq 2$ is a fixed integer, satisfy (Cf. Theorem
4.3): \be \sum_{k=1}^{n} \fr k{2l+k}(\laa_{n+2-k}-\laa_1)<
4(l-1)\laa_1. \en

In the final part of this paper, we will
prove universal inequalities for eigenvalues of
the polyharmonic operator of any order on compact domains with boundary in a unit sphere. For similar inequalities for eigenvalues of the Laplacian on compact domains in a sphere, we refer to [CY1], [AH] and the references therein.

\section{General Inequalities  for Eigenvalues of the Harmonic Operator of any Order on Riemannian Manifolds}

In this section, we  prove some general inequalities for eigenvalues of the polyharmonic operators
 on compact Riemannian manifolds.
\vskip0.2cm
{\bf Theorem 2.1.} {\it Let $(M,\ \langle\  ,\ \rangle )$ be an $n$-dimensional compact connected
Riemannian manifold with boundary $\pa M$ (possibly empty) and let $\nu$ be the outward unit normal vector
field of $\pa M$. Let $l$ be a positive integer and denote by $\Delta $ the Laplacian operator of $M$. Consider the eigenvalue problem
\setcounter{equation}{0}
\be
& &
(-\Delta)^l u=\la u\ \ \ {\rm in} \ \ M, \\
 \no
 & &  u|_{\pa M}=\left. \fr{\pa u}{\pa \nu}\right|_{\pa M}=\cdots =\left. \fr{\pa^{l-1} u}{\pa \nu^{l-1}}\right|_{\pa M}=0.
\en
Let  $\la_i,\ i=1,\cdots, $ be the $i$-th eigenvalue of the problem (2.1) and $u_i$ be the orthonormal eigenfunction corresponding to $\la_i$, that is,
\be
& & \no
 (-\Delta)^l u_i=\la_i u_i\ \ \ {\rm in} \  M, \\ \no
 & &
u_i|_{\pa M}=\left. \fr{\pa u_i}{\pa \nu}\right|_{\pa M}=\cdots= \left. \fr{\pa^{l-1} u_i}{\pa \nu^{l-1}}\right|_{\pa M}=0,\\ \no
& &
    \int_M u_i u_j=\delta_{ij}, \ \  {\rm for\ any\ \ } i,\ j=1, 2, \cdots.
\en
Then for any function $h\in C^{l+2}(M)\cap C^{l+1}(\pa M)$ and any positive integer $k$, we have
\be & &
\sum_{i=1}^k(\la_{k+1}-\la_i)^2\int_M hu_i\left((-\Delta)^l(hu_i)-\la_ihu_i\right)
\\ \no &\leq& \sum_{i=1}^k(\la_{k+1}-\la_i)||\left((-\Delta)^l(hu_i)-\la_i hu_i\right)||^2,
\en
\be
& &
\sum_{i=1}^k(\la_{k+1}-\la_i)^2\int_M\left(-hu_i^2\Delta h-2hu_i\langle\nabla h, \nabla u_i\rangle\right)\\ \no &\leq&  \delta\sum_{i=1}^{k}(\la_{k+1}-\la_i)^2\int_M hu_i\left((-\Delta)^l(hu_i)-\la_ihu_i)\right)+
\sum_{i=1}^k\fr{(\la_{k+1}-\la_i)}{\delta}\left|\left|\langle \nabla h, \nabla u_i\rangle +\fr{u_i\Delta h}2\right|\right|^2
\en
and
\be & &
\sum_{i=1}^k(\la_{k+1}-\la_i)^2\int_M\left(-hu_i^2\Delta h-2hu_i\langle\nabla h, \nabla u_i\rangle\right)\\ \no &\leq& \delta \sum_{i=1}^{k}(\la_{k+1}-\la_i)|| \left((-\Delta)^l(hu_i)-\la_ihu_i)\right) ||^2+\fr 1{\delta}
\sum_{i=1}^k(\la_{k+1}-\la_i)\left|\left|\langle \nabla h, \nabla u_i\rangle +\fr{u_i\Delta h}2\right|\right|^2,
\en
where
$\delta$ is any positive constant, $||g||^2=\int_M g^2$.
}
\vskip0.3cm
{\it Proof of Theorem 2.1.}  The inequality (2.2) follows from Theorem 2.1 in [AH]. In fact, by taking $N=1, \ B_1=h Id, \ A=(-\Delta)^l$
in Theorem 2.1 of [AH], one gets easily that
\be\no
\rho_i\equiv \langle [A, B_1]u_i, B_1u_i\ra=\int_M hu_i\left((-\Delta)^l(hu_i)-\la_ihu_i\right)
\en
and
\be\no
\Lambda_i\equiv|| [A, B_1]u_i||^2=||\left((-\Delta)^l(hu_i)-\la_i hu_i\right)||^2,
\en
which, by using (2.5) in [AH], gives (2.2).

We  use some  similar calculations as in [CY4] to prove (2.3). For $i=1,\cdots, k$,  consider the functions $\phi_i
: M\rightarrow I\!\!R$ given by
\be
\phi_i= hu_i-\sum_{j=1}^kr_{ij}u_j,
\en
where
\be
r_{ij}=\int_M hu_iu_j.
\en
Since
$$\phi_i|_{\pa M}=\left. \fr{\pa \phi_i}{\pa \nu}\right|_{\pa M}=\cdots =\left. \fr{\pa^{l-1} \phi_i}{\pa \nu^{l-1}}\right|_{\pa M}=0$$
and
\be\no
\int_M u_j\phi_{i}=0, \ \ \forall \ i,j=1,\cdots, k,
\en
it follows from the Rayleigh-Ritz  inequality that
\be
\la_{k+1}\int_M\phi_{i}^2&\leq &\int_M\phi_{i}(-\Delta)^l\phi_{i}\\ \no
&=&\la_i||\phi_i||^2+ \int_M\phi_{i}\left((-\Delta)^l\phi_{i}-\la_i hu_i\right)\\ \no
 &=&\la_i||\phi_i||^2+ \int_M\phi_{i}\left((-\Delta)^l(hu_i)-\la_i hu_i\right)\\ \no
 &=&\la_i||\phi_i||^2 +\int_M hu_{i}\left((-\Delta)^l(hu_i)-\la_i hu_i\right)-\sum_{j=1}^kr_{ij}s_{ij},
\en
where
\be\no
s_{ij}=\int_M \left((-\Delta)^l( hu_i)-\la_i hu_i\right)u_j.
\en
Notice that if $u\in C^{l+2}(M)\cap C^{l+1}(\pa M) $ satisfies
 \be
 u|_{\pa M}=\left. \fr{\pa u}{\pa \nu}\right|_{\pa M}=\cdots= \left. \fr{\pa^{l-1} u}{\pa \nu^{l-1}}\right|_{\pa M}=0,
\en
then
\be
 u|_{\pa M}&=&\left. \nabla u\right|_{\pa M}= \left. \Delta u\right|_{\pa M} =\left. \nabla(\Delta u)        \right|_{\pa M}  =\cdots =  \left. \Delta^{m-1} u\right|_{\pa M}\\ \no &=& \left. \nabla(\Delta^{m-1} u)\right|_{\pa M}=0 , \ \ \ \ \ \ \ \ \ \  {\rm when} \ \ \ l=2m
\en
and
\be
 u|_{\pa M}&=&\left. \nabla u\right|_{\pa M}= \left. \Delta u\right|_{\pa M} =\left. \nabla(\Delta u)        \right|_{\pa M}  =\cdots =  \left. \Delta^{m-1} u\right|_{\pa M}= \left. \nabla(\Delta^{m-1} u)\right|_{\pa M}\\ \no & = &\left. \Delta^{m} u\right|_{\pa M} =0 ,\ \ \ \ \ \ \ \ \ \  {\rm when} \ \ \ l=2m+1 .
\en
Observe that  both $u_j$ and $hu_i$ satisfy the boundary condition (2.8) and so they satisfy  (2.9) when $l=2m$ and  (2.10) when $l=2m+1$. Thus we can use
integration by parts to  conclude that
\be\no
\int_M u_j(-\Delta)^l(hu_i)=
\int_M hu_{i}(-\Delta)^l(u_j)=\la_j r_{ij},
\en
which gives
\be
s_{ij}=(\la_j-\la_i)r_{ij}.
\en
Set
\be
\no
p_i(h)=(-\Delta)^l(hu_i)-\la_ihu_i;
\en
then we have from (2.7) and (2.11) that
\be
(\la_{k+1}-\la_i)||\phi_{i}||^2
\leq \int_M \phi_ip_i(h)=
\int_M hu_ip_i(h)+\sum_{j=1}^k (\la_i-\la_j)r_{ij}^2.
\en
Set
\be
t_{ij}=\int_M u_j\left(\langle \nabla h, \nabla u_i\rangle+\fr{u_i \Delta h}2\right);
\en
then
$
t_{ij}+t_{ji}=  0$
and
\be \int_M (-2)\phi_i\left(\langle \nabla h, \nabla u_i\rangle+\fr{u_i \Delta h}2\right)=w_i+ 2\sum_{j=1}^k r_{ij}t_{ij},
\en
where
\be
w_i=\int_M (-h u_i^2\Delta h -2hu_i\langle\nabla h, \nabla u_i\rangle).
\en
Multiplying (2.14) by $(\la_{k+1}-\la_i)^2$ and using the Schwarz inequality and (2.12), we get
\be& &
(\la_{k+1}-\la_i)^2\left(w_i+ 2\sum_{j=1}^k r_{ij}t_{ij}\right)\\ \no &=&(\la_{k+1}-\la_i)^2
\int_M (-2)\phi_i\left(\left(\langle \nabla h, \nabla u_i\rangle+\fr{u_i \Delta h}2\right)-\sum_{j=1}^k t_{ij}u_j\right)
\\ \no &\leq&
\delta (\la_{k+1}-\la_i)^3||\phi_i||^2+\fr{(\la_{k+1}-\la_i)}{\delta}\int_M
\left|\langle \nabla h, \nabla u_i\rangle+\fr{u_i \Delta h}2-\sum_{j=1}^k t_{ij}u_j\right|^2
\\ \no &=&
\delta (\la_{k+1}-\la_i)^3||\phi_i||^2+\fr{(\la_{k+1}-\la_i)}{\delta}
\left(\left|\left|\langle \nabla h, \nabla u_i\rangle+\fr{u_i \Delta h}2\right|\right|^2-\sum_{j=1}^k t_{ij}^2\right)\\ \no &\leq&
\delta (\la_{k+1}-\la_i)^2\left(\int_M hu_ip_{i}(h)+\sum_{j=1}^k (\la_i-\la_j)r_{ij}^2\right)
\\ \no & & +\fr{(\la_{k+1}-\la_i)}{\delta}
\left(\left|\left|\langle \nabla h, \nabla u_i\rangle+\fr{u_i \Delta h}2\right|\right|^2-\sum_{j=1}^k t_{ij}^2\right).
\en
Summing over $i$  and noticing $r_{ij}=r_{ji}, \ t_{ij}=-t_{ji}$, we infer
\be \no
& &
\sum_{i=1}^k (\la_{k+1}-\la_i)^2w_i-2\sum_{i,j=1}^k (\la_{k+1}-\la_i) (\la_i-\la_j)r_{ij}t_{ij}
\\ \no &\leq&  \delta\sum_{i=1}^k (\la_{k+1}-\la_i)^2 \int_M hu_ip_{i}(h)
+\sum_{i=1}^k\fr{(\la_{k+1}-\la_i)}{\delta}
\left|\left|\langle \nabla h, \nabla u_i\rangle+\fr{u_i \Delta h}2\right|\right|^2
\\ \no & & - \sum_{i,j=1}^k (\la_{k+1}-\la_i)\delta(\la_i-\la_j)^2r_{ij}^2-\sum_{i,j=1}^k \fr{(\la_{k+1}-\la_i)}{\delta}t_{ij}^2.
\en Hence (2.3) is true. Substituting (2.2) into (2.3), one gets
(2.4).

\vskip0.3cm
We end this section by listing some Lemmas which  are needed in the next sections.
\vskip0.3cm
{\bf Lemma 2.1.} {\it Let $u_i$ and $\la_i, i=1, 2, \cdots , $ be as in Theorem 2.1, then
\be
0\leq \int_M u_i(-\Delta )^k u_i\leq \la_i^{k/l}, \ \ k=1,\cdots, l-1.
\en
}

{\it Proof.} When $k\in \{1,\cdots, l-1\}$ is even, we have
\be \no
\int_M u_i(-\D)^k u_i =\int_M u_i\D^k u_i=\int_M \left(\D^{k/2}u_i\right)^2\geq 0.
\en
On the other hand, if $k\in \{1,\cdots, l-1\}$ is odd,
\be \no
\int_M u_i(-\D)^k u_i &=& -\int_M u_i\D^k u_i
\\ \no
&=& -\int_M \D^{(k-1)/2}u_i\D\left(\D^{(k-1)/2}u_i\right)\\ \no
&=&
\int_M \left|\nabla\left( \D^{(k-1)/2}u_i\right)\right|^2
\\ \no &\geq& 0.
\en
Thus the inequality at the left hand side of (2.17) holds.

 We claim  that  for any $k=1,\cdots, l-1$,
 \be
 \left(\int_M u_i(-\Delta )^k u_i\right)^{k+1}\leq \left(\int_M u_i(-\Delta )^{k+1} u_i\right)^k.
\en
Since
\be \no
\left(\int_M u_i\Delta u_i\right)^2\leq \int_M u_i^2\int_M (\Delta u_i)^2=\int_M u_i \Delta^2 u_i,
\en
we know that (2.18) holds when $k=1$.

Suppose that (2.18) holds for $k-1$, that is
 \be
\left(\int_M u_i(-\Delta )^{k-1} u_i\right)^k\leq \left(\int_M u_i(-\Delta )^{k} u_i\right)^{k-1}.
\en
When $k$ is even, we have
\be
\int_M u_i(-\Delta )^k u_i&=&\int_M \Delta^{k/2-1}u_i \Delta\left(\Delta^{k/2}u_i\right)
\\ \no &= &  -\int_M \left\lan\nabla \left(\Delta^{k/2-1}u_i\right), \nabla\left(\Delta^{k/2}u_i\right)\right\ra
\\ \no
&\leq& \left(\int_M\left|\nabla \left(\Delta^{k/2-1}u_i\right)\right|^2\right)^{1/2}\left(\int_M\left|\nabla\left(\Delta^{k/2}u_i\right)\right|^2\right)^{1/2}\\ \no
&=& \left(-\int_M\Delta^{k/2-1}u_i\Delta^{k/2}u_i\right)^{1/2}\left(-\int_M\Delta^{k/2}u_i\Delta^{k/2+1}u_i\right)^{1/2}\\ \no
&=& \left(\int_M u_i(-\Delta)^{k-1}u_i\right)^{1/2}\left(\int_M u_i(-\Delta)^{k+1}u_i\right)^{1/2},
\en
On the other hand, when $k$ is odd,
\be
\int_M u_i(-\Delta )^k u_i&=&\int_M (-\Delta)^{(k-1)/2}u_i (-\Delta)^{(k+1)/2}u_i
\\ \no
&\leq& \left(\int_M \left((-\Delta)^{(k-1)/2}u_i\right)^2\right)^{1/2}\left(\int_M\left(\left(-\Delta \right)^{(k+1)/2}u_i\right)^2\right)^{1/2}\\ \no
&=& \left(\int_M u_i(-\Delta)^{k-1}u_i\right)^{1/2}\left(\int_M u_i(-\Delta)^{k+1}u_i\right)^{1/2}.
\en
Thus we always have
\be
\int_M u_i(-\Delta )^k u_i\leq  \left(\int_M u_i(-\Delta)^{k-1}u_i\right)^{1/2}\left(\int_M u_i(-\Delta)^{k+1}u_i\right)^{1/2}.
\en
Substituting (2.19) into (2.22), we know that (2.18) is true for $k$. Using (2.18) repeatedly, we get
\be
\int_M u_i(-\Delta )^k u_i\leq \left(\int_M u_i(-\Delta )^{k+1} u_i\right)^{k/(k+1)}
\leq \cdots \leq \left(\int_M u_i(-\Delta )^{l} u_i\right)^{k/l}=\la_i^{k/l}.
\en
This shows that the inequality at the right hand side of (2.17) also holds.

\vskip0.3cm {\bf Lemma 2.2.} {\it Let
$$C=\left\{z=(z_1,\cdots, z_n)\in {I\!\!R}^n|z_i\geq 0,\ i=1,\cdots, n,
\sum_{j=1}^n z_j=1\right\}.$$   Consider the function $f:
C\rightarrow R$ defined by \be f((z_1,\cdots, z_n))=
\sum_{i=1}^n\fr{z_i^2}{1+4z_i}.\en Then \be \min_{z\in C}
f(z)=f\left(\left(\fr 1n,\cdots, \fr 1n\right)\right)=\fr 1{n+4}.
\en}
\vskip0.3cm
{\it Proof of Lemma 2.2.} We minimize the function
$$\sum_{i=1}^n\fr{z_i^2}{1+4z_i}$$
with the constraint
$$\sum_{j=1}^n z_j=1,\  z_i\geq 0,\ i=1,\cdots, n. $$
By means of the method of the Lagrange multiplier, we consider the
following function:
$$ g=
\sum_{i=1}^n\fr{z_i^2}{1+4z_i}+\la \left(\sum_{j=1}^n z_j-1\right),
$$
where $\la $ is the Lagrange multiplier. The minimum point of
$\sum_{i=1}^n z_i^2/(1+4z_i)$ is a critical point of $g$. Taking the
derivative of $g$ with respect to $z_i$, we have
$$
\fr{2z_i(1+4z_i)-4z_i^2}{(1+4z_i)^2}+\la =0. $$ Multiplying the
above equaltion by $(1+4z_i)^2$ and simplifying, we get
$$
(16\la +4)z_i^2+(8\la +2)z_i+\la =0.
$$
Hence at most two of the $z_i's$ are distinct with each other at
critical point of $g$. Assume without loss of generality that
$z_1=z_2=\cdots =z_p=s,\ z_{p+1}=\cdots =z_{p+q}=t$ with $p+q=n$.
Then $ ps+qt=1$ and so we have \be  \sum_{i=1}^n
\fr{z_i^2}{(1+4z_i)}&=&\fr{ps^2}{1+4s}+\fr{qt^2}{1+4t}\\
\no &=& \fr{ps^2+qt^2+4st(ps+qt)}{(1+4s)(1+4t)}
\\
\no &=& \fr{ps^2+qt^2+4st}{(1+4s)(1+4t)}\\ \no &=&
\fr{(p+q+4)(ps^2+qt^2+4st(ps+qt))}{(n+4)(1+4s)(1+4t)}
\\ \no &=&\fr{p^2s^2+q^2t^2+pq(s^2+t^2)+4ps^2+4qt^2+4(p+q)st+16st}{(n+4)(1+4s)(1+4t)}
\\ \no &=&\fr{1-2pqst+pq(s^2+t^2)+4ps^2+4qt^2+4(p+q)st+16st}{(n+4)(1+4s)(1+4t)}
\\ \no &\geq&\fr{1+4ps^2+4qt^2+4(p+q)st+16st}{(n+4)(1+4s)(1+4t)}
\\ \no &=&\fr{1+4(ps+qt)(s+t)+16st}{(n+4)(1+4s)(1+4t)}
=\fr 1{n+4}. \en This completes the proof of Lemma 2.2.

\vskip0.3cm {\bf Lemma 2.3.}  Let $\{a_i\}_{i=1}^m$,
$\{b_i\}_{i=1}^m$ and $\{c_i\}_{i=1}^m$ be three  sequences of
non-negative real numbers with $\{a_i\}$ decreasing  and $\{b_i\}$
and $\{c_i\}_{i=1}^m$ increasing. Then the following inequality
holds: \be \left(\sum_{i=1}^m a_i^2b_i\right)\left(\sum_{i=1}^m a_i
c_i\right)\leq \left(\sum_{i=1}^m a_i^2 \right)\left(\sum_{i=1}^m
a_i b_i c_i\right). \en
{\it Proof.} When $m=1$, (2.27) holds
trivally. Suppose that (2.27) holds when  $m=k$, that is \be
\left(\sum_{i=1}^k a_i^2b_i\right)\left(\sum_{i=1}^k a_i
c_i\right)\leq \left(\sum_{i=1}^k a_i^2\right)\left(\sum_{i=1}^k a_i
b_i c_i\right). \en Then when $m=k+1$, we have from (2.28) that
\be
& & \left(\sum_{i=1}^{k+1} a_i^2\right) \left(\sum_{i=1}^{k+1} a_i
b_ic_i\right)-\left(\sum_{i=1}^{k+1} a_i^2b_i\right)\left(\sum_{i=1}^{k+1} a_i c_i\right) \\
\no &=& \left(\sum_{i=1}^{k} a_i^2\right) \left(\sum_{i=1}^k a_i b_i
c_i\right)-\left(\sum_{i=1}^{k} a_i^2b_i\right)\left(\sum_{i=1}^k
a_i c_i\right) + a_{k+1}^2\sum_{i=1}^k a_i b_i c_i
 \\ \no & & -a_{k+1}^2 b_{k+1}\sum_{i=1}^ka_i c_i+a_{k+1}b_{k+1}c_{k+1}\sum_{i=1}^k a_i^2-a_{k+1}c_{k+1}\sum_{i=1}^k a_i^2b_i\\ \no & \geq&  a_{k+1}^2\sum_{i=1}^k a_i b_i c_i
 -a_{k+1}^2 b_{k+1}\sum_{i=1}^ka_i c_i+a_{k+1}b_{k+1}c_{k+1}\sum_{i=1}^k a_i^2-a_{k+1}c_{k+1}\sum_{i=1}^k a_i^2b_i\\ \no &=&
 -a_{k+1}^2\sum_{i=1}^k(b_{k+1}-b_i)a_ic_i+a_{k+1}c_{k+1}\sum_{i=1}^k a_i^2(b_{k+1}-b_i)\\ \no
 &=&\sum_{i=1}^k a_{k+1}a_i(b_{k+1}-b_i)(c_{k+1}a_i-a_{k+1}c_i)
 \\ \no &\geq&0.
 \en
 Where in the last inequality we have used the fact that
 \be\no
  a_{k+1}a_i(b_{k+1}-b_i)(c_{k+1}a_i-a_{k+1}c_i)\geq 0, \ \ i=1,\cdots,k.
 \en
 Thus (2.27) holds for $m=k+1$. This completes the proof of Lemma 2.3.
\vs
The following result is the so called {\it Reverse Chebyshev Inequality} (Cf. [HLP]).
\vskip0.3cm
{\bf Lemma 2.4.} Suppose
$\{a_i\}_{i=1}^m$ and $\{b_i\}_{i=1}^m$ are two real sequences with
$\{a_i\}$ increasing  and $\{b_i\}$ decreasing. Then the following
inequality holds: \be \sum_{i=1}^m a_ib_i\leq \fr 1m
\left(\sum_{i=1}^m a_i\right)\left(\sum_{i=1}^m b_i\right). \en
\vs
The following lemma can be also found in [HLP]
\vskip0.3cm {\bf Lemma 2.5.} {\it Let
$\{c_k\}_{k=1}^l$ and $\{d_k\}_{k=1}^l$ be two increasing real
sequences. Then for any permutation $\{i_1,\cdots, i_l\}$ of
$\{1,\cdots, l\}$, we have \be\no \sum_{k=1}^l c_k d_{i_k}\geq
c_1d_l+c_2d_{l-1}+\cdots+c_l d_1. \en}
\vskip0.3cm
{\bf Remark.} Lemma 2.4 also admits a probabilistic interpretation. We may assume that the $a_i$ and $b_i$ are nonnegative and satisfy $\sum_{i=1}^m a_i=1=\sum_{i=1}^m b_i$ so that we can interprete them as the probabilities for observing $i$ under the laws $a$ or $b$, resp. One then needs to prove $\sum_{i=1}^m a_ib_i\leq \fr 1m$. When the $b_i$ are all the same, that is, $=\fr 1m$, the inequality is obviously an equality, and when $b_i$ is decreasing instead of constant, the right hand side stays the same, but the left hand side can only become smaller, because then  higher weights are placed on those $i$ with smaller $a_i$. Thus, the inequality follows. In fact, Lemma 2.5 above admits a similar interpretation.

\section{ Universal Inequalities for Eigenvalues of the Polyharmonic Operators   on Compact Domains in $I\!\!R^n$}
\setcounter{equation}{0}
In this section, we will prove  universal bounds on eigenvalues of the polyharmonic operator on  bounded domains in a Euclidean space by using Theorem 2.1.
\vskip0.3cm
{\bf Theorem 3.1.} { \it Let $\om $ be a connected dounded domain in an  $n$-dimensional Euclidean space $I\!\!R^n$ and let  $\Delta$ be the Laplacian of $I\!\!R^n$. Denote by
$\la_i$ the $i$-th eigenvelue of the eigenvalue problem:
\be \no
& & (-\Delta)^l u=\la u\ \ \ {\rm in} \ \ \om, \\
\no
& &   u|_{\pa \om}=\left. \fr{\pa u}{\pa \nu}\right|_{\pa \om}=\cdots =\left. \fr{\pa^{l-1} u}{\pa \nu^{l-1}}\right|_{\pa \om}=0.
\en
Then we have
\be
& &
\sum_{i=1}^k\left(\la_{k+1}-\la_i\right)^2
\\ \no &\leq& \left(\frac{4l(n+2l-2)}{n^2}\right)^{1/2}
\left(\sum_{i=1}^k\left(\la_{k+1}-\la_i\right)^2\la_i^{(l-1)/l}\right)^{1/2} \left(\sum_{i=1}^k\left(\la_{k+1}-\la_i\right)\la_i^{1/l}\right)^{1/2}
\en
}
\vskip0.3cm
{\bf Corollary 3.1.} {\it Under the same assumptions as in Theorem 3.1, we have
\be
& & \la_{k+1}\\ \no
&\leq &\fr 1k\sum_{i=1}^k\la_i+\fr{2l(n+2l-2)}{k^2n^2}\left(\sum_{i=1}^k\la_i^{(l-1)/l}
\right)\left(\sum_{i=1}^k\la_i^{1/l}\right)\\ \no
& & +\left\{\left(\fr{2l(n+2l-2)}{k^2n^2}\right)^2\left(\sum_{i=1}^k\la_i^{(l-1)/l}
\right)^2\left(\sum_{i=1}^k\la_i^{1/l}\right)^2-\fr 1k\sum_{i=1}^k\left(\la_i-\fr 1k\sum_{j=1}^k \la_j\right)^2\right\}^{\fr 12}.
\en
and
\be
& & \la_{k+1}\\ \no
&\leq & \left(1+\fr{2l(n+2l-2)}{n^2}\right)\fr 1k\sum_{i=1}^k\la_i\\
\no & & +
\left\{\left(\fr{2l(n+2l-2)}{n^2}\fr 1k\sum_{i=1}^k\la_i\right)^2-\left(1+\fr{4l(n+2l-2)}{n^2}\right)\fr 1k\sum_{i=1}^k\left(\la_i-\fr 1k\sum_{j=1}^k\la_j\right)^2\right\}^{1/2}
\en
}
\vskip0.3cm
{\it Proof of Theorem 3.1. }
Let $x_1, x_2, \cdots, x_n$ be the standard Euclidean coordinate functions of $I\!\!R^n$. Let $u_i$  be the $i$-th orthonormal eigenfunction corresponding to the eigenvalue $\la_i$ of the problem (3.1), $i=1,\cdots$; then
\be
\Delta x_{\alpha}=0,\   \nabla x_{\alpha}=(\underbrace{0, \cdots, 0,  1}_{\alpha}, 0,\cdots, 0), \ \ \alpha =1,2,\cdots, n,
\en
which implies that
\be
(-\Delta )^l(x_{\alpha} u_i)&=&x_{\alpha} (-\Delta )^lu_i+2l (-1)^l\langle \nabla x_{\alpha}, \nabla(\Delta^{l-1}u_i)\rangle\\ \no
&=& \la_i x_{\alpha}u_i+2l (-1)^l\langle \nabla x_{\alpha}, \nabla(\Delta^{l-1}u_i)\rangle
\en
Taking $h=x_{\alpha}$ in (2.3), we infer for any $\delta>0$ that
\be & & \no
\sum_{i=1}^k(\la_{k+1}-\la_i)^2\int_{\om}\left(-2x_{\alpha}u_i\langle\nabla x_{\alpha}, \nabla u_i\rangle\right)\\ \no &\leq&
\delta\sum_{i=1}^k\left(\la_{k+1}-\la_i\right)^2\int_{\om}2l(-1)^lx_{\alpha}u_i\langle\na x_{\alpha}, \na(\D^{l-1}u_i)\ra \\ \no & &
+\fr1{\delta}\sum_{i=1}^k\left(\la_{k+1}-\la_i\right)||\langle\nabla x_{\alpha}, \nabla u_i\rangle||^2.
\en
Summing over $\alpha$, we have
\be & &
\sum_{i=1}^k(\la_{k+1}-\la_i)^2\sum_{\alpha=1}^{n}\int_{\om}\left(-2x_{\alpha}u_i\langle\nabla x_{\alpha}, \nabla u_i\rangle\right)\\ \no &\leq&
2l\delta\sum_{i=1}^k\left(\la_{k+1}-\la_i\right)^2\sum_{\alpha=1}^{n}\int_{\om}(-1)^lx_{\alpha}u_i \langle\na x_{\alpha}, \na(\D^{l-1}u_i)\ra
\\ \no & & +\fr1{\delta}\sum_{i=1}^k\left(\la_{k+1}-\la_i\right)\sum_{\alpha=1}^{n}||\langle\nabla x_{\alpha}, \nabla u_i\rangle||^2.
\en
Observe that
\be
\sum_{{\alpha}=1}^n|\nabla x_{\alpha}|^2 =n, \ \ \   \sum_{{\alpha}=1}^n\langle\nabla x_{\alpha},\nabla u_i\rangle^2=|\nabla u_i|^2.
\en
Hence
\be
\sum_{\alpha=1}^n\int_{\om}\left(-2x_{\alpha}u_i\langle\nabla x_{\alpha}, \nabla u_i\rangle\right)
=\fr12\sum_{\alpha=1}^n\int_{\om}u_i^2\Delta x_{\alpha}^2=\int_{\om}u_i^2\sum_{\alpha=1}^n|\nabla x_{\alpha}|^2=n,
\en
From (2.17), we infer
\be
\sum_{\alpha=1}^{n}||\langle\nabla x_{\alpha}, \nabla u_i\rangle||^2=\int_{\om}|\na u_i|^2=\int_{\om}(-u_i\D u_i)\leq \la_i^{1/l}.
\en
Since
\be\no
\D^{l-1}(x_{\alpha}u_i)=2(l-1)\langle\na x_{\alpha}, \na(\D^{l-2}u_i)\ra +x_{\alpha}\D^{l-1} u_i,
\en
we have
\be
\int_{\om}x_{\alpha}u_i \langle\na x_{\alpha}, \na(\D^{l-1}u_i)\ra  &=&\int_{\om}x_{\alpha}u_i \D^{l-1} \langle\na x_{\alpha}, \na u_i\ra\\ \no &=&\int_{\om}\D^{l-1}(x_{\alpha}u_i)  \langle\na x_{\alpha}, \na u_i\ra\\ \no &=&\int_{\om}\left(2(l-1)\langle\na x_{\alpha}, \na(\D^{l-2}u_i)\ra +x_{\alpha}\D^{l-1} u_i\right)\lan \na x_{\alpha}, \na u_i\ra.
\en
On the other hand,
\be
\int_{\om}x_{\alpha}u_i \langle\na x_{\alpha}, \na(\D^{l-1}u_i)\ra &=&  -\int_{\om}\D^{l-1}u_i\ {\rm div} (x_{\alpha} u_i\na x_{\alpha})\\ \no &=& -\int_{\om}\D^{l-1}u_i(|\na x_{\alpha}|^2u_i+x_{\alpha}\lan \na x_{\alpha}, \na u_i\ra),
\en
where ${\rm div} (X)$ denotes the divergence of $X$.
Combining  (3.10) and (3.11), we obtain
\be & &
\int_{\om}x_{\alpha}u_i \langle\na x_{\alpha}, \na(\D^{l-1}u_i)\ra
 \\ \no &=&
 \int_M\left\{(l-1)\langle\na x_{\alpha}, \na (\D^{l-2}u_i)\ra \lan \na x_{\alpha}, \na u_i\ra-\fr 12\D^{l-1}u_i |\na x_{\alpha}|^2u_i\right\}
\en
It then follows from (2.17), (3.7), (3.12) and
\be\no
\sum_{\alpha=1}^{n}\langle\na x_{\alpha}, \na (\D^{l-2}u_i)\ra\langle\na x_{\alpha}, \na u_i\ra=\langle\na u_i, \na (\D^{l-2}u_i)\ra
\en
that
\be
& &
\sum_{\alpha=1}^n\int_{\om}(-1)^lx_{\alpha}u_i \langle\na x_{\alpha}, \na(\D^{l-1}u_i)\ra\\ \no &=& \sum_{\alpha=1}^n\int_{\om}(-1)^l\left\{(l-1)\lan \na x_{\alpha}, \na(\D^{l-2}u_i)\ra\lan \na x_{\alpha}, \na u_i\ra -\frac 12 \D^{l-1}u_i|\na x_{\alpha}|^2 u_i\right\}
\\ \no &=& \int_{\om}(-1)^l\left\{(l-1)\lan \na(\D^{l-2}u_i), \na u_i\ra -\frac n2 u_i\D^{l-1}u_i\right\}
\\ \no &=&
\left(l-1+\fr n2\right)\int_{\om} u_i(-\D)^{l-1}u_i \\ \no &\leq &
\left(l-1+\fr n2\right)\la_i^{(l-1)/l}.
\en
Substituting (3.8), (3.9) and (3.13) into (3.6), one gets
\be
& & n\sum_{i=1}^k\left(\la_{k+1}-\la_i\right)^2
\\ \no & \leq& l(n+2l-2)\delta \sum_{i=1}^k\left(\la_{k+1}-\la_i\right)^2\la_i^{(l-1)/l} +\fr 1{\delta}\sum_{i=1}^k\left(\la_{k+1}-\la_i\right)\la_i^{1/l}.
\en
Taking
\be \no
\delta=\left\{\fr{\sum_{i=1}^k\left(\la_{k+1}-\la_i\right)\la_i^{1/l}}{
l(n+2l-2) \sum_{i=1}^k\left(\la_{k+1}-\la_i\right)^2\la_i^{(l-1)/l}}\right\}^{1/2},
\en
we get (3.1).
\vs
In the proof of Corollary 3.1 we will use the reverse-Chebyshev inequality which was used earlier for similar purposes in [A1] and [AH].
\vskip0.3cm
{\it Proof of Corollary 3.1.}
It follows from (2.30)  that
\be
\sum_{i=1}^k\left(\la_{k+1}-\la_i\right)\la_i^{1/l}\leq\fr 1k\left(\sum_{i=1}^k\left(\la_{k+1}-\la_i\right)\right)\left(\sum_{i=1}^k\la_i^{1/l}\right)
\en
and
\be
\sum_{i=1}^k\left(\la_{k+1}-\la_i\right)^2\la_i^{(l-1)/l}\leq \fr 1k\left(
\sum_{i=1}^k\left(\la_{k+1}-\la_i\right)^2\right)\left(\sum_{i=1}^k\la_i^{(l-1)/l}\right).
\en
Introducing (3.15) and (3.16) into (3.1), we infer
\be\no
& &
\sum_{i=1}^k\left(\la_{k+1}-\la_i\right)^2
\\ \no &\leq& \frac{4l(n+2l-2)}{k^2n^2}
\left(\sum_{i=1}^k\left(\la_{k+1}-\la_i\right)\right)\left(\sum_{i=1}^k\la_i^{(l-1)/l}\right) \left(\sum_{i=1}^k\la_i^{1/l}\right),
\en
Solving this quadratic polynomial about $\la_{k+1}$, one gets (3.2).

From (2.27), we have \be \no& & \left(
\sum_{i=1}^k\left(\la_{k+1}-\la_i\right)^2\la_i^{(l-1)/l}\right)\left(\sum_{i=1}^k\left(\la_{k+1}-\la_i\right)\la_i^{1/l}\right)
\\ \no & \leq& \left(\sum_{i=1}^k\left(\la_{k+1}-\la_i\right)^2\right)\left(\sum_{i=1}^k\left(\la_{k+1}-\la_i\right)\la_i\right).
\en It then follows from (3.1) that \be
\sum_{i=1}^k\left(\la_{k+1}-\la_i\right)^2\leq \fr{4l(n+2l-2)}{n^2}
\sum_{i=1}^k\left(\la_{k+1}-\la_i\right)\la_i, \en which implies
(3.3).

\section{ Universal Inequalities for Lower Order Eigenvalues of the Polyharmonic Operators on Compact Domains in  $I\!\!R^n$}
\setcounter{equation}{0}
 In [AB4], Ashbaugh and Benguria showed that when
$l=1$, the first $n+1$ eigenvalues of the problem (1.1) satisfies
the inequality $\la_2+\la_3+\cdots +\la_{n+1}\leq (n+4)\la_1$. Also,
Ashbaugh showed in [A1] that when $l=2$, $\la_2+\la_3+\cdots
+\la_{n+1}\leq (n+24)\la_1$. The following result generalizes the
estimate by Ashbaugh and Benguria to any $l$ and strengthens the
above Ashbaugh's inequality.
\vskip0.3cm{\bf Theorem 4.1.} {\it
Under the same assumptions as in Theorem 3.1, we have \be
\sum_{i=2}^{n+1}\la_i+\sum_{i=1}^{n-1}\fr{2(l-1)i}{2l+i-1}(\la_{n+1-i}-\la_1)\leq(n+4l(2l-1))\la_1.
\en } \vskip0.3cm {\it Proof of Theorem 3.2.}
As in the proof of Theorem 3.1, we let $u_i$  be the $i$-th orthonormal eigenfunction corresponding to the eigenvalue $\la_i$ of the problem (3.1), $i=1,\cdots$.
We first {\it claim}  that
there exists a set of Cartesian coordinate system  $(x_1,...,x_n)$ of $I\!\!R^n$
so that the following orthogonality conditions are satisfied: \be
\int_{\om}x_iu_1u_j=0\ \ \ {\rm for}\ \ \ 1\leq j\leq i\leq n. \en
 Indeed, by choosing the
origin properly, we can assume that there exists a Cartesian
coordinates $(y_1,...,y_n)$ of $R^n$ such that \be
\int_{\om}y_iu_1^2=0\ \ \ {\rm for}\ \ \ i=1,...,n.\en Consider the
matrix $A$ defined by
$$
A=\left[\begin{array}{cccc} \int_{\om}y_1u_1u_2 &
\int_{\om}y_1u_1u_3 & \cdots &\int_{\om}y_1u_1u_{n+1}\\
\int_{\om}y_2u_1u_2 & \int_{\om}y_2u_1u_2 & \cdots
&\int_{\om}y_2u_1u_{n+1}\\ \cdots & \cdots & \cdots & \cdots\\
\int_{\om}y_nu_1u_2 & \int_{\om}y_nu_1u_3 & \cdots
&\int_{\om}y_nu_1u_{n+1}
\end{array}\right]
$$
From the orthogonalization of Gram-Schmidt(QR-factorization
theorem), we know that $A$ can be written as
$$
B=TA, $$ where $T=(t_{ij})$ is an orthogonal $n\times n$ matrix and
$B$ is an upper triangular matrix. Hence, we have, for any $k$ and
$j$ with $k>j$,
$$
\sum_{l=1}^n  t_{kl}\int_{\om} y_l u_1 u_{j+1}=0. $$ Defining new
coordinate  functions $x_k,$ by $x_k=\sum_{j=1}^n t_{kj}y_j$, one
has, for any $i, j=1,2,...,n,$ satisfying $i>j$, \be \int_{\om}
x_iu_1 u_{j+1}=0. \en Combining (4.3) and (4.4), we know that our
{\it claim} is true.

Since (4.2) holds, for each $i=1,...,n,$ we get  from the
Rayleigh-Ritz inequality, (3.5) and (3.12) that \be
\la_{i+1}&\leq&\fr{\int_{\om}x_iu_1(-\Delta)^l(x_i u_1)}{\int_{\om}
x_i^2 u_1^2}\\ \no &=& \la_1+\fr{ 2l\int_{\om} (-1)^l x_iu_1\langle
\na x_i, \na(\Delta^{l-1}u_1)\rangle}{\int_{\om} x_i^2 u_1^2}\\
\no &=& \la_1+\fr{ 2l\int_{\om} (-1)^l\left\{(l-1)\langle\na x_i,
\na (\D^{l-2}u_1)\ra \lan \na x_i, \na u_1\ra-\fr 12 u_1\D^{l-1}u_1
\right\}}{\int_{\om} x_i^2 u_1^2}. \en
Since
$$1=\int_{\om}u_1^2=-2\int_{\om}x_iu_1\lan\na x_i, \na u_1\ra,$$
we have by squaring both sides and using the Cauchy-Schwarz
inequality that \be \no 1\leq 4\int_{\om}\lan \na x_i, \na
u_1\ra^2\int_{\om} x_i^2 u_1^2, \en which gives \be \fr 1{\int_{\om}
x_i^2 u_1^2}\leq 4\int_{\om}\lan \na x_i, \na u_1\ra^2. \en
Substituting (4.6) into (4.5) yields \be& & \la_{i+1}-\la_1\\ \no
&\leq &\left(2l\int_{\om} (-1)^l\left\{(l-1)\langle\na x_i, \na
(\D^{l-2}u_1)\ra \lan \na x_i, \na u_1\ra-\fr 12 u_1\D^{l-1}u_1
\right\}\right)\left(4\int_{\om}\lan \na x_i, \na u_1\ra^2\right).
\en Set \be a=\int_{\om} u_1(-\D)^{l-1}u_1 , \ \ a_i=\int_{\om}
(-1)^l\langle\na x_i, \na(\D^{l-2} u_1)\ra \lan \na x_i, \na
u_1\rangle, \ \ i=1,\cdots, n; \en then \be \sum_{i=1}^na_i=a.\en
Take a permutation $\{i_1,\cdots, i_n\}$ of $\{1,\cdots, n\}$ so
that \be a_{i_n}\leq a_{i_{n-1}}\leq\cdots \leq a_{i_1}, \en it then
follows from (4.9) and (4.10) that \be a_{i_k}\leq \fr 1k a, \
k=1,\cdots, n.\en Substituting (4.11) into (4.7), we get \be
\la_{i_k+1}-\la_1\leq 8l\left(\fr{(l-1)}k a +\fr a2
\right)\int_{\om}\lan \na x_{i_k}, \na u_1\ra^2. \en Multiplying
(4.12) by $\fr{(2l-1)k}{2(l-1)+k}$ and simplifying, one has \be
\left(1+\fr{2(l-1)(k-1)}{2(l-1)+k}\right)(\la_{i_k+1}-\la_1)&=&
\fr{(2l-1)k}{2(l-1)+k}( \la_{i_k+1}-\la_1)\\ \no &\leq&
4l(2l-1)a\int_{\om}\lan \na x_{i_k}, \na u_1\ra^2. \en Summing over
$k$,
we have \be & & 4l(2l-1)a\int_{\om}|\na u_1|^2\\ \no &\geq&
\sum_{k=1}^n\left(1+\fr{2(l-1)(k-1)}{2(l-1)+k}\right)(\la_{i_k+1}-\la_1)\\
\no &=&\sum_{k=1}^n(\la_k-\la_1)+\sum_{k=2}^n
\fr{2(l-1)(k-1)}{2(l-1)+k}\left(\la_{i_k+1}-\la_1\right). \en
Observe that $\{i_1,\cdots, i_n\}$ is a permutation of $\{1,\cdots,
n\}$. We {\it claim } that there is a permutation $\{q_2,
q_3,\cdots, q_{n}\}$ of $\{1,...,n-1\}$ such that \be \sum_{k=2}^{n}
\fr{2(l-1)(k-1)}{2(l-1)+k}\left(\la_{i_k+1}-\la_1\right) \geq
\sum_{k=2}^{n}
\fr{2(l-1)(k-1)}{2(l-1)+k}\left(\la_{q_k+1}-\la_1\right). \en In
fact, if $i_1=n$, then $\{i_2,\cdots, i_n\}$ is a permutation of
$\{1,\cdots, n-1\}$ and there is nothing to prove. On the other
hand, if $i_1=m\in \{1, 2,\cdots, n-1\}$, then $\{i_2, i_3,\cdots,
i_{n}\}=\{1, 2,\cdots, m-1, m+1,\cdots, n\}$ and so there is a $j\in
\{2,\cdots, n\}$ such that $i_j=n$ which implies that \be\{i_2,
i_3,\cdots, i_{j-1}, i_{j+1},\cdots, i_{n}\}=\{1, 2,\cdots, m-1,
m+1,\cdots, n-1\}, \en and \be & & \sum_{k=2}^{n}
\fr{2(l-1)(k-1)}{2(l-1)+k}\left(\la_{i_k+1}-\la_1\right)\\ \no
&=&\sum_{k=2, k\neq j}^{n}
\fr{2(l-1)(k-1)}{2(l-1)+k}\left(\la_{i_k+1}-\la_1\right)
+\fr{2(l-1)(j-1)}{2(l-1)+j}\left(\la_{n+1}-\la_1\right)\\ \no &\geq&
\sum_{k=2, k\neq
j}^{n}\fr{2(l-1)(k-1)}{2(l-1)+k}\left(\la_{i_k+1}-\la_1\right)
+\fr{2(l-1)(j-1)}{2(l-1)+j}\left(\la_{m+1}-\la_1\right). \en From
(4.16), we know that $\{i_2, i_3,\cdots, i_{j-1}, m, i_{j+1},\cdots,
i_{n}\}$ is a permutation of $\{1, 2,\cdots, m-1, m,  m+1,\cdots,
n-1\}.$ Set $i_2=q_2,\cdots,i_{j-1}=q_{j-1}, m=q_j,
i_{j+1}=q_{j+1},\cdots, i_{n}=q_{n}$; then $\{q_2, q_3,\cdots,
q_{n}\}$ is a permutation of $\{1,...,n-1\}$ and  we can rewrite
(4.17) as (4.15). Thus our {\it claim} is true.

Since $\left\{\fr{2(l-1)(k-1)}{2(l-1)+k}\right\}_{k=2}^n$ and
$\{\la_k-\la_1\}_{k=2}^n$ are two increasing sequences and $\{
q_1+1, q_2+1,\cdots, q_n+1\}$ is a permutation of $\{2,\cdots, n\}$,
we conclude from Lemma 2.5 that \be \sum_{k=2}^{n}
\fr{2(l-1)(k-1)}{2(l-1)+k}\left(\la_{q_k+1}-\la_1\right)
\geq\sum_{k=2}^{n} \fr{2(l-1)(k-1)}{2(l-1)+k}(\la_{n-k+2}-\la_1).
\en
 Thus we
have from (2.17),  (4.14) and (4.18) that \be & &
\sum_{k=1}^n\la_{i+1}+\sum_{k=2}^n\fr{2(l-1)(k-1)}{2(l-1)+k}(\la_{n-k+2}-\la_1)
\\ \no &\leq& n\la_1+4l(2l-1) a\int_{\om}|\na u_1|^2 \\ \no &=&
n\la_1+4l(2l-1)\left(\int_{\om}
u_1(-\D)^{l-1}u_1\right)\left(\int_{\om}u_1(-\D u_1)\right)\\
\no &\leq&
n\la_1+4l(2l-1)\la_1^{(l-1)/l}\cdot\la_1^{1/l}=(n+4l(2l-1))\la_1.
\en This is just the inequality (4.1). The  proof of Theorem 3.2 is
completed.
\vskip0.3cm
Our next result is to prove the inequality (1.16) as mentioned in the introduction.
\vskip0.3cm {\bf Theorem 4.2.}
{\it Let $\om$ be a
connected bounded domain with smooth boundary in  ${I\!\!R}^n$. Denote by $\nu$ the outward unit normal vector field of $\pa \om$ and let
$\laa_i, \ i=1,\cdots, n+1,$ be the first $(n+1)$ eigenvalues of the following buckling
problem: \be \Delta^2 u=-\laa \Delta u \ \ \ {\rm in} \ \ \om ,  \ \
\ \ \left. u\right|_{\pa \om}=\left. \fr{\pa u}{\pa \nu}\right|_{\pa
\om}=0. \en Then, \be \sum_{i=1}^n\laa_{i+1}
+\fr{4(\laa_2-\laa_1)}{n+4}\leq (n+4)\laa_1. \en }
\vskip0.3cm {\it Proof of Theorem 4.2.} Let us denote by  $u_i$  be the $i$-th orthonormal eigenfunction corresponding to the eigenvalue $\la_i$ of the buckling problem (4.20), $i=1,\cdots$. That is, we have
 \be \Delta^2 u_i=-\laa_i \Delta u \ \ \ {\rm in} \ \ \om ,  \ \
\ \ \left. u_i\right|_{\pa \om}=\left. \fr{\pa u_i}{\pa \nu}\right|_{\pa
\om}=0, \\
\int_{\om}\lan\na u_i,  \na u_j\ra=\delta_{ij}, \ \ \forall  \ \ i, j=1,\cdots.
\en
Using similar  discussions as in the proof of Theorem 4.1, we can
find a set of Cartesian coordinates $(x_1,...,x_n)$ of $I\!\!R^n$
so
that the following orthogonality conditions are satisfied:
\be
\int_{\om}\lan\na(x_iu_1), \na u_j\ra=-\int_{\om} x_i u_1 \Delta u_j=0\ \ \ {\rm for}\ \ \ 1\leq j\leq i\leq n.
\en
Now we start with the well known Rayleigh-Ritz inequality
\be
\laa_{i+1}\leq \fr{\int_{\om} \phi \D^2 \phi}{\int_{\om}|\na
\phi|^2}, \en which is satisfied by any sufficiently smooth function
$\phi $ such that \be\no & & \phi =\fr{\partial \phi}{\partial
\nu}=0\ \ {\rm on\ \ } \pa \om,\\  & & \no\int_{\om}\lan\na \phi, \na
u_j\ra=0, \ j=1,\cdots, i. \en
Set $u=u_1$ and we choose as our trial function $$ \phi =x_i u, $$
which clearly satisfies the above boundary condition, and by (4.24)
also the orthogonality condition. Thus we have \be \laa_{i+1}
\int_{\om}|\na (x_i u)|^2\leq \int_{\om}x_iu\D^2(x_i u),\ \
i=1,\cdots, n. \en Let us calculate \be \int_{\om}|\na (x_i
u)|^2=\int_{\om} x_i^2|\na u|^2+ 2\int_{\om}x_i u u_{x_i}+\int_{\om}
u^2=\int_{\om} x_i^2 |\na u|^2, \en
where $u_{x_i}=\langle \na x_{i}, \na u\rangle=\fr{\pa u}{\pa x_{i}}$.
As for the right hand side of
(4.26), we have \be \int_{\om} x_i u\D^2(x_i u)&=& \int_{\om} x_i
u\left(x_i \D^2 u+ 4\D u_{x_i}\right)\\ \no &=&
-\laa_1\int_{\om}x_i^2u \D u+ 4\int_{\om} x_i u\D u_{x_i}. \en
By integration by parts, one gets \be \int_{\om}x_i^2u \D u&=&
-\int_{\om}\lan\na(x_i^2 u), \na u\ra\\ \no &=& -\int_{\om} x_i^2|\na
u|^2-2\int_{\om} x_i u u_{x_i}\\
\no &=& -\int_{\om} x_i^2|\na u|^2 +\int_{\om} u^2. \en Thus, \be
\int_{\om} x_i u\D^2(x_i u)&=& \laa_1\int_{\om}x_i^2|\na u|^2
-\laa_1\int_{\om} u^2+ 4 \int_{\om} x_i u \D u_{x_i}. \en
Substituting (4.27) and (4.30) into (4.26) and dividing both sides by
$\int_{\om}x_i^2|\na u|^2$, we get \be
\laa_{i+1}-\laa_1\leq\fr{-\laa_1\int_{\om} u^2+ 4 \int_{\om} x_i u\D
u_{x_i}}{\int_{\om}x_i^2|\na u|^2}.\en We have
 \be
\int_{\om} x_i u \D u_{x_i}&=&-\int_{\om}\lan\na(x_i u), \na u_{x_i}\ra\\
\no &=& -\int_{\om} (x_i \lan\na u,  \na u_{x_i}\ra + u u_{x_i x_i})\\ \no
&=& \fr 12\int_{\om}|\na u|^2+\int_{\om} u_{x_i}^2\\ \no &=& \fr 12
+\int_{\om} u_{x_i}^2. \en The Cauchy-Schwarz inequality implies \be
1&=&\int_{\om}(-u\D u)= \left(\int_{\om}(-u\D u)\right)^2 \leq
\left(\int_{\om} u^2\right)\left(\int_{\om}(\D u)^2\right)
\\ \no &=& \left(\int_{\om}
u^2\right)\left(\int_{\om}u \D^2 u\right)=\left(\int_{\om}
u^2\right)\laa_1 \int_{\om}|\na u|^2=\laa_1 \int_{\om} u^2.
 \en
Substituting (4.32) and (4.33) into (4.31) , we get \be
\laa_{i+1}-\laa_1\leq\fr{1+ 4
\int_{\om}u_{x_i}^2}{\int_{\om}x_i^2|\na u|^2}.\en
Since
\be
\left(\int_{\om}\lan\na(x_i u), \na u_{x_i}\ra\right)^2 &\leq&
\left(\int_{\om}|\na (x_i u)|^2\right)\left(\int_{\om} |\na
u_{x_i}|^2\right)\\ \no &=&\left(\int_{\om}x_i^2|\na
u|^2\right)\left(\int_{\om} |\na u_{x_i}|^2\right), \en it follows
that \be\fr{\left(\int_{\om}\lan\na(x_i u), \na u_{x_i}\ra\right)^2
}{\int_{\om}x_i^2|\na u|^2}\leq \int_{\om} |\na u_{x_i}|^2.\en
 Combining (4.32) and (4.36), one gets
\be\fr{1+4\int_{\om}u_{x_i}^2+4\left(\int_{\om}u_{x_i}^2\right)^2
}{\int_{\om}x_i^2|\na u|^2}\leq 4\int_{\om} |\na u_{x_i}|^2.\en Set
\be b_i=\int_{\om}u_{x_i}^2 {\ \ \ \rm   and\ \ \ }
\epsilon_i=\fr{\left(\int_{\om}u_{x_i}^2\right)^2}{1+4\int_{\om}u_{x_i}^2}.
\en
 It then follows from (4.34) and (4.37) that \be
(1+4\epsilon_i)(\laa_{i+1}-\laa_1)\leq 4\int_{\om} |\na
u_{x_i}|^2.\en
Since
$u|_{\pa\om}=\left.\fr{\pa u}{\pa \nu}\right|_{\pa \om}=0$, we know that $u_{x_i}|_{\pa\om}=0$, which implies from the divergence theorem that
\be\no
\int_{\om}|\na u_{x_i}|^2=-\int_{\om}u_{x_i}\D u_{x_i}=
-\int_{\om}u_{x_i}(\D u)_{x_i}.
\en
Let
$e_i=(\underbrace{0, \cdots, 0,  1}_{i}, 0,\cdots, 0)$ and consider the vector field $X=u(\Delta u)_{x_i} e_i$. We infer from the divergence theorem and $X|_{\pa \om}=0$ that
$$0=\int_{\om} {\rm div} X=\int_{\om} (u_{x_i}(\Delta u)_{x_i}+u(\Delta u)_{x_ix_i}).
$$
Hence, we have
\be \int_{\om}
|\na u_{x_i}|^2=\int_{\om} u(\D
u)_{x_ix_i}. \en
Substituting (4.40) into (4.39) and summing on $i$ from
1 to $n$, we have
 \be
\sum_{i=1}^n(1+4\epsilon_i)(\laa_{i+1}-\laa_1)&\leq& 4\int_{\om}
u\left(\sum_{i=1}^n(\D u)_{x_ix_i}\right)\\ \no
&=&-4\laa_1\int_{\om}u\D u\\ \no
&=& 4 \laa_1\int_{\om}|\na u|^2=4\laa_1.\en Now we want to estimate
the left hand side of the above inequality. We have
 \be
\sum_{i=1}^n(1+4\epsilon_i)(\laa_{i+1}-\laa_1)&\geq &
\sum_{i=1}^n(\laa_{i+1}-\laa_1)+ 4(\laa_{2}-\laa_1) \sum_{i=1}^n
\epsilon_i. \en

From the definition, we know that \be \sum_{i=1}^n b_i=\int_{\om}
|\na u|^2=1.\en
Thus we deduce from Lemma 2.2 that \be
\sum_{i=1}^n\epsilon_i =f((b_1,\cdots, b_n))\geq \fr 1{n+4}. \en
Combining (4.41), (4.42) and (4.44), we get \be \sum_{i=2}^{n+1}\laa_i
+\fr{4(\laa_2-\laa_1)}{n+4}\leq (n+4)\laa_1. \en This completes the
proof of Theorem 4.2.
\vskip0.3cm
The final result of this section is to prove the inequality (1.18). That is, we have
\vskip0.3cm {\bf Theorem 4.3.} {\it Let $l\geq 2$ be a
positive integer and let $\om$ be a connected bounded domain with
smooth boundary in  ${I\!\!R}^n$.
Consider the eigenvalue problem
 \be
& & (-\Delta)^l u=-\laa\Delta u\ \ \ {\rm in} \ \ \om, \\
\no  & &  u|_{\pa \om}=\left. \fr{\pa u}{\pa \nu}\right|_{\pa
\om}=\cdots =\left. \fr{\pa^{l-1} u}{\pa \nu^{l-1}}\right|_{\pa
\om}=0. \en Let \be\no
 0<\laa_1\leq\laa_2\leq\cdots\leq \laa_{n+1},
\en denote the first $n+1$ eigenvalues of the above problem. Then we
have \be \sum_{k=1}^{n} \fr k{2l+k}(\laa_{n+2-k}-\laa_1)<
4(l-1)\laa_1. \en} \vskip0.3cm {\it Proof of Theorem 4.3.} Let $v_i$
be the $i$-th orthonormal eigenfunction of the problem (4.46)
corresponding to the eigenvalue $\laa_i$, $i=1, 2, \cdots,$ that is,
$v_i$ satisfies \be\no & & (-\Delta)^lv_i=-\laa_i \Delta v_i \ \
{\rm in \ \ } \om,
\\ \no & &
\ \left. v_i\right|_{\pa \om}=\left. \fr{\pa v_i}{\pa
\nu}\right|_{\pa \om}=0=\cdots =\left. \fr{\pa^{l-1} v_i}{\pa
\nu^{l-1}}\right|_{\pa \om}=0, \\ \no & &  \int_{\om}\langle \nabla
v_i, \nabla v_j\rangle=\delta_{ij},\ \ \ \forall\ i, j. \en
As in the proof of Theorem 4.1, we take a
set of Cartesian coordinates $(x_1,...,x_n)$ of $I\!\!R^n$
so that the
following orthogonality conditions are satisfied: \be
\int_{\om}\lan\na(x_iv_1), \na v_j\ra= - \int_{\om}x_iv_1\D v_j=0\ \ \ {\rm
for}\ \ \ 1\leq j\leq i\leq n. \en
The Rayleigh-Ritz inequality now
states that
\be \laa_{i+1}\leq \fr{\int_{\om} \phi (-\D)^l
\phi}{\int_{\om}|\na \phi|^2}, \en which is satisfied by any
sufficiently smooth function $\phi $ such that \be \no & & \phi
=\fr{\partial \phi}{\partial \nu}=\cdots =\fr{\partial^{l-1}
\phi}{\partial \nu^{l-1}}=
0\ \ {\rm on\ \ } \pa \om,\\
\no & & \int_{\om}\lan\na \phi, \na v_j\ra=0, \ j=1,\cdots, i. \en

Set $v=v_1$ and we choose as our trial function $ \phi =x_i v, $
which clearly satisfies the above boundary condition, and by (4.48)
also the orthogonality condition. Then we have \be \laa_{i+1}\leq
\fr{\int_{\om} x_iv (-\D)^l(x_i v) }{\int_{\om}|\na(x_i v) |^2},\
i=1,\cdots, n. \en Since $\Delta x_{i}=0, i=1,2,\cdots, n, $ we have
 \be (-\Delta )^l(x_{i}
v)&=&x_{i} (-\Delta )^lv+2l
(-1)^l (\Delta^{l-1}v)_{x_i}\\
\no &=& \laa_1 x_{i}(-\D v)+2l (-1)^l (\Delta^{l-1}v)_{x_i}. \en As
calculated in (4.27) and (4.29),  we have \be
\int_{\om}|\na(x_iv)|^2=\int_{\om}x_i^2 |\na v|^2, \
\int_{\om}x_i^2v \D v= -\int_{\om} x_i^2|\na v|^2 +\int_{\om} v^2.
\en
Substituting (4.51) and (4.52) into (4.50), we get
\be
\laa_{i+1}-\laa_1\leq\fr{-\laa_1\int_{\om} v^2+ 2l (-1)^l
\int_{\om}x_iv(\Delta^{l-1}v)_{x_i}}{\int_{\om}x_i^2|\na v|^2}. \en
Since \be\no \D^{l-1}(x_{i}v)=2(l-1)
 (\D^{l-2}v)_{x_i}
+x_{i}\D^{l-1}v, \en we have \be \int_{\om}x_{i}v
(\Delta^{l-1}v)_{x_i} &=&\int_{\om}x_iv \D^{l-1}v_{x_i}\\
\no &=&\int_{\om}\D^{l-1}(x_iv) v_{x_i}\\ \no &=&\int_{\om}(2(l-1)
 (\D^{l-2}v)_{x_i}+x_{i}\D^{l-1}v)v_{x_i}\\ \no &=&\int_{\om}(-2(l-1)
 \D^{l-2}v v_{x_ix_i}+x_{i}\D^{l-1}v v_{x_i}).\en
On the other hand, we have from the divergence theorem that
\be\int_{\om}x_{i}v
(\Delta^{l-1}v)_{x_i} = -\int_{\om}\D^{l-1}v( v+x_iv_{x_i}). \en Combining (4.54) and
(4.55), we obtain \be\int_{\om}x_{i}v
(\Delta^{l-1}v)_{x_i}&=&\int_{\om}\left((l-1)
 (\D^{l-2}v)_{x_i}v_{x_i}-\fr 12 v\D^{l-1}v\right)\\ \no &=&\int_{\om}\left(-(l-1)
 (\D^{l-2}v)_{x_ix_i}v-\fr 12 v\D^{l-1}v\right).
\en
 Substituting (4.56) into (4.53), we infer
 \be
\laa_{i+1}-\laa_1\leq\fr{-\laa_1\int_{\om} v^2+ l
\int_{\om}(-2(l-1)((-\Delta)^{l-2}v)_{x_ix_i}v+v(-\D)^{l-1}v)}{\int_{\om}x_i^2|\na
v|^2}. \en By using the same arguments as in the proof of (4.37),
one deduces \be\no
\fr{1+4\int_{\om}v_{x_i}^2+4\left(\int_{\om}v_{x_i}^2\right)^2
}{\int_{\om}x_i^2|\na v|^2}\leq 4\int_{\om} |\na v_{x_i}|^2, \en
which gives \be \fr{1}{\int_{\om}x_i^2|\na v|^2}\leq4\int_{\om} |\na
v_{x_i}|^2 \en and since
$\sum_{i=1}^n\int_{\om}v_{x_i}^2=\int_{\om}|\na v|^2=1$, we know
that there exists at least one  $i\in\{1,\cdots, n\}$ such that
(4.58) is a strict inequality. Set \be a=\int_{\om} v(-\D)^{l-1}v, \
\ a_i=\int_{\om}(-((-\Delta)^{l-2}v)_{x_ix_i}v), \ i=1,\cdots n;\en
then \be \sum_{i=1}^na_i=a.\en Introducing (4.59) into (4.57), we
have \be \laa_{i+1}-\laa_1\leq 4\left(-\laa_1\int_{\om} v^2+
l(2(l-1)a_i+a)\right) \int_{\om} |\na v_{x_i}|^2. \en and there
exists at least $i\in\{1,\cdots, n\}$ such that (4.61) is a strict
inequality. Take a permutation $\{i_1,\cdots, i_n\}$ of $\{1,\cdots,
n\}$ so that \be a_{i_n}\leq a_{i_{n-1}}\leq\cdots \leq a_{i_1}, \en
it then follows from (4.60) and (4.62) that \be a_{i_k}\leq \fr 1k
a, \ k=1,\cdots, n.\en Substituting (4.63) into (4.61), we get \be
\laa_{i_k+1}-\laa_1&\leq& 4\left(-\laa_1\int_{\om} v^2+
l(2(l-1)a_{i_k}+a)\right) \int_{\om} |\na v_{x_{i_k}}|^2\\ \no
&\leq& 4\left(-\laa_1\int_{\om} v^2+
l\left(\fr{2(l-1)}k+1\right)a\right) \int_{\om} |\na v_{x_{i_k}}|^2
\en and for some $k\in \{1,\cdots, n\}$, (4.64) is a strict
inequality. Before we can finish the proof of Theorem 4.3, let us
prove the following inequalities :
 \be & & \int_{\om}
v(-\Delta )^{l-1} v\leq \laa_1^{(l-2)/(l-1)}, \\
 & & \int_{\om} v\D^2v\leq \laa_1^{1/(l-1)}, \\ & &  \int_{\om}
v^2\geq \laa_1^{-1/(l-1)}. \en
First observe as in the proof of Lemma 2.1 that for any $k=1,\cdots, l-1, $
 \be\no
 \int_{\om} v(-\D)^kv\geq 0.
 \en
 When $l=2$, (4.65) and (4.66) hold obviously and in this case we have from Schwarz inequality that
\be \no 1&=&\int_{\om}(-v\D v)= \left(\int_{\om}(-v\D v)\right)^2
\leq \left(\int_{\om} v^2\right)\left(\int_{\om}(\D v)^2\right)=
\laa_1\int_{\om} v^2.\en
Hence (4.67) holds when $l=2$.

Assume now that $l>2.$ We
claim  that  for any $k=2,\cdots, l-1$,
 \be \left(\int_{\om} v(-\Delta
)^k v\right)^{k}\leq \left(\int_{\om} v(-\Delta )^{k+1}
v\right)^{k-1}. \en Since \be \no \int_{\om} v\Delta^2v =\int_{\om}
\Delta v\Delta v =-\int_{\om} \na \Delta v\na v,\en
 we have from
Schwarz inequality that
 \be  \left(\int_{\om} v\Delta^2 v\right)^2
 \leq
 \left(\int_{\om} |\na \D v|^2\right)\left(\int_{\om} |\na v|^2\right)
= -\int_{\om} \D v\D^2 v = \int_{\om} v(-\D^3 v). \en
Hence (4.68) holds when $k=2$. Suppose that
(4.68) holds for $k-1$, that is
 \be
\left(\int_{\om} v(-\Delta )^{k-1} v\right)^{k-1}\leq
\left(\int_{\om} v(-\Delta  )^{k} v\right)^{k-2}. \en
As in the proof of (2.26), we  have \be
\int_{\om} v(-\Delta )^k v\leq  \left(\int_{\om}
v(-\Delta)^{k-1}v\right)^{1/2}\left(\int_{\om}
v(-\Delta)^{k+1}v\right)^{1/2}. \en Substituting (4.70) into (4.71),
we know that (4.68) is true for $k$. Using (4.68) repeatedly, we get
\be \no \int_{\om} v(-\Delta )^k v\leq \left(\int_{\om} v(-\Delta
)^{k+1}v\right)^{(k-1)/k} \leq \cdots \leq \left(\int_{\om}
v(-\Delta )^{l} v\right)^{(k-1)/(l-1)}=\laa_1^{(k-1)/(l-1)}. \en
Taking $k=2$ and $k=l-1$ in the above inequality, respectively, one gets (4.65) and (4.66). On the other
hand, we have from Schwarz inequality and (4.66) that
\be 1&=&\int_{\om}(-v\D v)= \left(\int_{\om}(-v\D v)\right)^2
\leq \left(\int_{\om} v^2\right)\left(\int_{\om}(\D v)^2\right)
\\ \no &=& \left(\int_{\om}
v^2\right)\left(\int_{\om}v \D^2 v\right)\leq
\laa_1^{1/(l-1)}\left(\int_{\om} v^2\right).\en This proves (4.67).
Now we continue on the proof of Theorem 4.3. Substituting (4.65) and
(4.67) into (4.64) and multiplying both sides by $\fr k{k+2l}$, we
get \be \fr k{k+2l}(\laa_{i_k+1}-\laa_1) \leq 4(l-1)
\laa_1^{(l-2)/(l-1)}\int_{\om} |\na v_{x_{i_k}}|^2 \en and for some
$k\in\{1,\cdots, n\}$, the above inequality is a strict inequality.
Thus by summing on $k$ and using $\int_{\om} v\D^2 v\leq
\laa_1^{1/(l-1)}$, one gets
 \be
\sum_{k=1}^n \fr k{2l+k}(\laa_{i_k+1}-\laa_1) &<& 4(l-1)
\laa_1^{(l-2)/(l-1)}\sum_{k=1}^n \int_{\om} |\na v_{x_{i_k}}|^2\\
\no &=&
4(l-1) \laa_1^{(l-2)/(l-1)}\sum_{k=1}^n \int_{\om} |\na v_{x_{k}}|^2\\
\no &=&  4(l-1) \laa_1^{(l-2)/(l-1)} \int_{\om} v\D^2v\leq
4(l-1)\laa_1. \en
 Since $\left\{\fr k{2l+k}\right\}_{k=1}^{n}$ and
$\{\laa_{k+1}-\laa_1\}_{k=1}^{n}$ are two increasing sequences, we
have from Lemma 2.5 that \be \sum_{k=1}^{n} \fr
k{2l+k}(\laa_{i_k+1}-\laa_1)\geq
 \sum_{k=1}^{n} \fr k{2l+k}(\laa_{n+2-k}-\laa_1).\en
Substituting (4.75)  into (4.74), we get (4.47). This completes the
proof of Theorem 4.3.

\section{ Eigenvalues of the Polyharmonic Operators on Compact  Domains in a Unit Sphere}
\setcounter{equation}{0}
In this section, we will prove  universal inequalities for eigenvalues of the polyharmonic operators on compact  connected domains in a unit $n$-sphere $S^n$. Let $l$ be a positive integer and for $p=0, 1, 2,...,$  define the polynomials $F_p(t)$ inductively by
\be
& & F_0(t)=1,\ \ \ F_1(t)=t-n,\ \ \ \\ \no
& & F_p(t)=(2t-2)F_{p-1}(t)-(t^2+2t-n(n-2))F_{p-2}(t), \ \ p=2,\cdots.
\en
Set
\be
F_l(t)=t^l+a_{l-1}t^{l-1}+\cdots +a_1t+a_0.
\en
\vskip0.3cm
{\bf Theorem 5.1.} { \it Let $\la_i$ be the i-th eigenvalue of the following eigenvalue problem:
\be \no
& & (-\Delta)^l u=\la u \ \ \ {\rm in} \ \ \om, \\ \no
& &
 \left. u\right|_{\pa \om}=\left. \fr{\pa u}{\pa \nu}\right|_{\pa \om}=\cdots=\left.\fr{\pa^{l-1} u}{\pa \nu^{l-1}}\right|_{\pa \om}=0,
\en
where $\om$ is a compact connected domain in a unit $n$-sphere $S^n$.
Then we have
\be& &
\sum_{i=1}^k\left(\la_{k+1}-\la_i\right)^2\\ \no &\leq&
\fr 1n\left\{
\sum_{i=1}^k(\la_{k+1}-\la_i)^2\left(|a_{l-1}|\la_i^{(l-1)/l}+\cdots +|a_1|\la_i^{1/l}+|a_0|\right)\right\}^{1/2}
\\ \no
& & \times \left\{
\sum_{i=1}^k(\la_{k+1}-\la_i)\left(n^2+4\la_i^{1/2}\right)\right\}^{1/2}.
\en
}
\vskip0.3cm
{\bf Corollary 5.1.} {\it Under the same assumptions as in Theorem 5.1, we have
\be
& & \la_{k+1}
\\ \no &\leq& \fr 1k\sum_{i=1}^k \la_i+\fr 1{2n^2k^2}\left(\sum_{i=1}^k\left(|a_{l-1}|\la_i^{(l-1)/l}+\cdots +|a_1|\la_i^{1/l}+|a_0|\right)\right)\left(kn^2+4\sum_{i=1}^k\la_i^{1/l}\right)\\ \no
& & +\left\{\fr 1{4n^4k^4}\left(\sum_{i=1}^k\left(|a_{l-1}|\la_i^{(l-1)l}+\cdots +|a_1|\la_i^{1/l}+|a_0|\right)\right)^2\left(kn^2+4\sum_{i=1}^k\la_i^{1/l}\right)^2\right. \\ \no & & \ \ \ \left. -\fr 1k\sum_{i=1}^k\left(\la_i-\fr 1k\sum_{j=1}^k\la_j\right)^2\right\}^{1/2}
\en
and
 \be
\la_{k+1}\leq  U_{k+1}+\sqrt{U_{k+1}^2-V_{k+1}},
\en
where
\be
U_{k+1}=\fr 1k\sum_{i=1}^k \la_i+\fr 1{2n^2k}\sum_{i=1}^k\left(|a_{l-1}|\la_i^{(l-1)/l}+\cdots +|a_1|\la_i^{1/l}+|a_0|\right)
\left(n^2+4\la_i^{1/l}\right)
\en
and
\be V_{k+1}=\fr 1k\sum_{i=1}^k\la_i^2+\fr 1{n^2k}\sum_{i=1}^k\la_i\left(|a_{l-1}|\la_i^{(l-1)/l}+\cdots +|a_1|\la_i^{1/l}+|a_0|\right)
\left(n^2+4\la_i^{1/l}\right).
\en
}

{\it Proof of Theorem 5.1.}   Let $x_1, x_2,\cdots, x_{n+1}$ be the standard  coordinate functions of the Euclidean space  $I\!\!R^{n+1}$; then
$$S^n=\{(x_1,\dots, x_{n+1})\in  I\!\!R^{n+1}; \sum_{\alpha=1}^{n+1}x_{\alpha}^2=1\}.
$$ It is well known that
\be
\Delta x_{\alpha}=-nx_{\alpha},\ \ \  \alpha =1,\cdots, n+1.
\en

Let $u_i$ be the $i$-th orthonormal eigenfunction  corresponding to the eigenvalue $\la_i$, $i=1, 2, \cdots.$
For any $\delta> 0$, by taking $h=x_{\alpha}$ in (2.4),  we have
\be& & \no
\sum_{i=1}^k(\la_{k+1}-\la_i)^2\int_{\om}(-x_{\alpha}u_i^2\Delta x_{\alpha}-2x_{\alpha}u_i\langle\nabla x_{\alpha}, \nabla u_i\rangle)
\\ \no &\leq& \delta\sum_{i=1}^{k+1}(\la_{k+1}-\la_i)^2\int_{\om} x_{\alpha}u_i
((-\Delta )^l(x_{\alpha}u_i)-\la_i x_{\alpha}u_i)\\ \no & & +\fr
1{\delta}\sum_{i=1}^k(\la_{k+1}-\la_i)\left|\left|\langle \nabla
x_{\alpha}, \nabla u_i\rangle +\fr{u_i\Delta
x_{\alpha}}2\right|\right|^2 \en Taking sum on $\alpha$ from $1$ to
$n+1$, we get \be& & \sum_{i=1}^k(\la_{k+1}-\la_i)^2\sum_{\alpha
=1}^{n+1}\int_{\om}\left(-x_{\alpha}u_i^2\Delta
x_{\alpha}-2x_{\alpha}u_i\langle\nabla x_{\alpha}, \nabla
u_i\rangle\right)\\ \no &\leq&
\delta\sum_{i=1}^{k+1}(\la_{k+1}-\la_i)^2\sum_{\alpha
=1}^{n+1}\int_{\om} x_{\alpha}u_i ((-\Delta )^l(x_{\alpha}u_i)-\la_i
x_{\alpha}u_i)\\ \no & & +\fr
1{\delta}\sum_{i=1}^k(\la_{k+1}-\la_i)\sum_{\alpha
=1}^{n+1}\left|\left|\langle \nabla x_{\alpha}, \nabla u_i\rangle
+\fr{u_i\Delta x_{\alpha}}2\right|\right|^2. \en Using
$\sum_{\alpha=1}^{n+1}x_{\alpha}^2=1$, (2.17) and (5.8), we infer
\be & & \sum_{\alpha =1}^{n+1}\int_{\om}\left(-x_{\alpha}u_i^2\Delta
x_{\alpha}-2x_{\alpha}u_i\langle\nabla x_{\alpha}, \nabla
u_i\rangle\right)
\\ \no &=& \int_{\om}\left(\left(\sum_{\alpha =1}^{n+1}x_{\alpha}^2\right)nu_i^2
-u_i\left\langle\nabla\left(\sum_{\alpha =1}^{n+1}x_{\alpha}^2\right), \nabla u_i\right\rangle\right)=\int_{\om}nu_i^2=n.
\en
\be & &
\sum_{\alpha=1}^{n+1}\left|\left|\langle \nabla x_{\alpha}, \nabla u_i\rangle +\fr{u_i\Delta x_{\alpha}}2\right|\right|^2
\\ \no &=& \int_{\om}\sum_{\alpha=1}^{n+1}\left(\langle \nabla x_{\alpha}, \nabla u_i\rangle^2-n\langle \nabla x_{\alpha}, \nabla u_i\rangle u_ix_{\alpha}+\fr{n^2u_i^2x_{\alpha}^2}4\right)\\ \no
&= &
\fr{n^2}4+\int_{\om}|\nabla u_i|^2\\ \no & =&
\fr{n^2}4+\int_{\om}u_i(-\Delta u_i)
\\ \no
& \leq& \fr{n^2}4+\la_i^{1/l}.
\en

For any smooth function $f$ on $\om$, we have from the Bochner formula that
\be
\fr 12\Delta|\nabla f|^2 &=& |\nabla^2f|^2+\langle \nabla f, \nabla(\Delta f)\rangle + {\rm Ric}(\nabla f, \nabla f)\\ \no  &=& |\nabla^2f|^2+\langle \nabla f, \nabla(\Delta f)\rangle + (n-1)|\nabla f|^2,
\en
where ${\rm Ric} $ is the Ricci tensor of $S^n$. Thus for any smooth function $g$ on $\om$,
\be
\fr 12\Delta|\nabla g|^2 = |\nabla^2g|^2+\langle \nabla g, \nabla(\Delta g)\rangle +
(n-1)|\nabla g|^2
\en
and
\be
\fr 12\Delta|\nabla (f+g)|^2 = |\nabla^2(f+g)|^2+\langle \nabla (f+g), \nabla(\Delta (f+g))\rangle +
(n-1)|\nabla (f+g)|^2.
\en
Subtracting the sum of (5.12) and (5.13) from (5.14), we get
\be
 \Delta \langle \nabla f, \nabla g\rangle =2\langle \nabla^2 f, \ \nabla^2 g\rangle
 +\langle \nabla f, \nabla(\Delta g)\rangle +\langle \nabla g, \nabla(\Delta f)\rangle
+2(n-1)\langle\nabla f, \ \nabla g\rangle,
\en
where
\be
\no
\langle \nabla^2 f, \ \nabla^2 g\rangle =\sum_{s, t=1}^n\nabla^2 f(e_s, e_t)\nabla^2 g(e_s, e_t),
\en
being $e_1, \cdots, e_n$ orthonormal vector fields locally defined on $\om$.
Since
\be \no
\nabla^2 x_{\alpha}=-x_{\alpha} \langle\ ,  \ \rangle,
\en
we infer from (5.15) by taking $f=x_{\alpha}$ that
\be
\Delta \langle \nabla x_{\alpha}, \nabla g\rangle&=&-2x_{\alpha}\Delta g+\langle \nabla x_{\alpha}, \nabla(\Delta g)\rangle+(n-2)\langle\nabla x_{\alpha}, \ \nabla g\rangle
\\ \no &=& -2x_{\alpha}\Delta g+\langle \nabla x_{\alpha}, \nabla((\Delta +(n-2))g)\rangle.
\en
For each $q=0, 1,\cdots$, thanks to (5.8) and (5.16), there are polynomials $B_q$ and $C_q$ of degrees less than or equal to $q$ such that
\be
\D^q(x_{\alpha} g)=x_{\alpha}B_q(\D)g+2\langle \nabla x_{\alpha}, \nabla(C_q(\D)g)\rangle.
\en
It is obvious that
\be
B_0=1, \ \ B_1=t-n,\ \ C_0=0,\ \ C_1=1.
\en
It follows from (5.8), (5.16)  and (5.17) that
\be& &
\D^q(x_{\alpha} g)\\ \no
&=&\D(\D^{q-1}(x_{\alpha}g))\\
\no
&=&\D(x_{\alpha}B_{q-1}(\D)g+2\langle \nabla x_{\alpha}, \nabla(C_{q-1}(\D)g)\rangle)\\ \no
&=&x_{\alpha}((\D-n)B_{q-1}(\D)-4\D C_{q-1}(\D))g+2\langle\nabla x_{\alpha}, \nabla((B_{q-1}(\D)+(\D+(n-2))C_{q-1}(\D))g)\rangle.
\en
Thus, for any $q=2,\cdots$,  we have
\be& &
B_{q}(\D)=(\D-n)B_{q-1}(\D)-4\D C_{q-1}(\D), \\
& &
C_{q}(\D)=B_{q-1}(\D)+(\D+(n-2))C_{q-1}(\D).
\en
Consequently, we have
\be\ \ \ \
B_q(\D)&=&(2\D-2)B_{q-1}(\D)-(\D+n-2)B_{q-1}(\D)-4\D C_{q-1}(\D)\\ \no
&=&(2\D-2)B_{q-1}(\D)-(\D+n-2)((\D-n)B_{q-2}(\D)-4\D C_{q-2}(\D))-4\D C_{q-1}(\D)\\ \no
&=&(2\D-2)B_{q-1}(\D)-(\D^2+2\D-n(n-2))B_{q-2}(\D)\\
\no & & +4\D[B_{q-2}(\D)+(\D+n-2)C_{q-2}(\D)-C_{q-1}(\D)]
\\ \no
&=&(2\D-2)B_{q-1}(\D)-(\D^2+2\D-n(n-2))B_{q-2}(\D), \ \ q=2,\cdots .
\en
Since (5.18) and (5.22) hold, we know that $B_q=F_q, \ \ \forall q=0, 1,\cdots .$
It follows from (5.17) and the divergence theorem that
\be& &
\int_{\Omega} x_{\alpha}u_i((-\D)^l(x_{\alpha}u_i)-\la_ix_{\alpha}u_i)
\\ \no &=&\int_{\Omega} x_{\alpha}u_i\left((-1)^l\left( x_{\alpha}B_l(\D)u_i+2\langle \nabla x_{\alpha},
\nabla(C_l(\D)u_i)\rangle\right)-\la_i x_{\alpha} u_i\right)
\\ \no &=&\int_{\Omega} x_{\alpha}u_i\left((-1)^l\left(x_{\alpha}(\D^l+a_{l-1}\D^{l-1}+\cdots+a_0)u_i+2\langle\nabla x_{\alpha}, \nabla(C_l(\D)u_i)\rangle\right)-\la_ix_{\alpha}u_i\right)
\\ \no &=&\int_{\Omega} (-1)^lx_{\alpha}u_i\left(x_{\alpha}(a_{l-1}\D^{l-1}+\cdots+a_0)u_i+2\langle\nabla x_{\alpha}, \nabla(C_l(\D)u_i)\rangle\right)
\en
Summing  on $\alpha$, one has from $\nabla\left(\sum_{\alpha=1}^{n+1}x_{\alpha}^2\right)=0$ and (2.17) that
\be& & \int_{\Omega} x_{\alpha}u_i((-\D)^l(x_{\alpha}u_i)-\la_ix_{\alpha}u_i)
\\ \no
&=&\int_{\Omega} u_i(-1)^l(a_{l-1}\D^{l-1}+\cdots+a_0)u_i
\\ \no
&\leq& |a_{l-1}|\int_{\Omega} u_i(-\D)^{l-1}u_i+\cdots +|a_1|\int_{\Omega} u_i(-\D)u_i+
|a_0|\int_{\Omega} u_i^2\\ \no
&\leq& |a_{l-1}|\la_i^{(l-1)/l}+ \cdots +|a_1|\la_i^{1/l}+|a_0|.
\en
Substituting (5.10), (5.11) and (5.24) into (5.9), we infer
\be\no
n\sum_{i=1}^k(\la_{k+1}-\la_i)^2&\leq& \delta \sum_{i=1}^k(\la_{k+1}-\la_i)^2(|a_{l-1}|\la_i^{(l-1)/l}+ \cdots +|a_1|\la_i^{1/l}+|a_0|)\\ \no
& & +\fr 1{\delta}\sum_{i=1}^k(\la_{k+1}-\la_i)
\left(\la_i^{1/l}+\fr{n^2}4\right).
\en
Taking
\be
\delta =\left\{\fr{\sum_{i=1}^k(\la_{k+1}-\la_i)
\left(\la_i^{1/l}+\fr{n^2}4\right)}{\sum_{i=1}^k(\la_{k+1}-\la_i)^2(|a_{l-1}|\la_i^{(l-1)/l}+ \cdots +|a_1|\la_i^{1/l}+|a_0|)}\right\}^{1/2},
\en
we get  (5.1). This completes the proof of Theorem 5.1.
\vskip0.3cm
{\it Proof of Corollary 5.1.}
From Lemma 2.4, we have
\be & &
\sum_{i=1}^k(\la_{k+1}-\la_i)^2\left(|a_{l-1}|\la_i^{(l-1)/l}+ \cdots +|a_1|\la_i^{1/l}+|a_0|\right)
\\ \no &\leq& \fr 1k \left\{
\sum_{i=1}^k(\la_{k+1}-\la_i)^2\right\}\left\{
\sum_{i=1}^k\left(|a_{l-1}|\la_i^{(l-1)/l}+ \cdots +|a_1|\la_i^{1/l}+|a_0|\right)\right\}
\en
and
\be
\sum_{i=1}^k(\la_{k+1}-\la_i)\left(n^2+4\la_i^{1/l}\right)
&\leq&\fr 1k\left\{
\sum_{i=1}^k(\la_{k+1}-\la_i)\right\}
\left\{
\sum_{i=1}^k\left(n^2+4\la_i^{1/l}\right)\right\}\\ \no
&=& \fr 1k\left\{
\sum_{i=1}^k(\la_{k+1}-\la_i)\right\}
\left\{
kn^2+4\sum_{i=1}^k\la_i^{1/l}\right\}.
\en
Substituting (5.26) and (5.27) into (5.3), we get
\be & &
\sum_{i=1}^k(\la_{k+1}-\la_i)^2\\ \no
&\leq& \fr 1{n^2k^2}\left\{\sum_{i=1}^k
\left(|a_{l-1}|\la_i^{(l-1)/l}+ \cdots +|a_1|\la_i^{1/l}+|a_0|\right)\right\}
\left\{ \sum_{i=1}^k(\la_{k+1}-\la_i)\right\}\left\{
kn^2+4\sum_{i=1}^k\la_i^{1/l}\right\}.
\en
Solving this quadratic polynomial of $\la_{k+1}$, we have (5.4).

On the other hand, one gets by using (2.27) that
\be& &
\left\{
\sum_{i=1}^k(\la_{k+1}-\la_i)^2\left(|a_{l-1}|\la_i^{(l-1)/l}+\cdots +|a_1|\la_i^{1/l}+|a_0|\right)\right\}
\left\{
\sum_{i=1}^k(\la_{k+1}-\la_i)\left(n^2+4\la_i^{1/2}\right)\right\}
\\ \no &\leq&
\left\{\sum_{i=1}^k(\la_{k+1}-\la_i)^2\right\}\left\{\sum_{i=1}^k(\la_{k+1}-\la_i)
\left(|a_{l-1}|\la_i^{(l-1)/l}+\cdots +|a_1|\la_i^{1/l}+|a_0|\right)
\left(n^2+4\la_i^{1/2}\right)\right\}, \en which, combining with
(5.3), gives \be& & \sum_{i=1}^k(\la_{k+1}-\la_i)^2\\ \no & &
\leq\fr 1{n^2}
\sum_{i=1}^k(\la_{k+1}-\la_i)\left(|a_{l-1}|\la_i^{(l-1)/l}+ \cdots
+|a_1|\la_i^{1/l}+|a_0|\right)\left(n^2+4\la_i^{1/l} \right). \en
Hence \be \la_{k+1}\leq U_{k+1}+\sqrt{U_{k+1}^2-V_{k+1}}, \en where
$U_{k+1}$ and $V_{k+1}$ are given by (5.6) and (5.7), respectively.
Thus (5.5) holds.
\vskip0.3cm {\it Acknowledgements.}  The fourth
author  thanks the Max Planck Institute for Mathematics
in the Sciences for its hospitality during the preparation of this
paper and CAPES.

\noindent J\"urgen Jost ( jost@mis.mpg.de ), Xianqing Li-Jost (
xli-jost@mis.mpg.de )

\noindent Max Planck Institute for Mathematics in the Sciences,
04103 Leipzig, Germany

\vskip0.3cm \noindent Qiaoling Wang (
wang@mat.unb.br )

\noindent Departamento de Matem\'atica, UnB, 70910-900,
Bras\'{\i}lia-DF, Brazil

\vskip0.3cm
\noindent Changyu Xia ( xia@mat.unb.br ) \noindent
Departamento de Matem\'atica, UnB, 70910-900, Bras\'{\i}lia-DF,
Brazil

\noindent and Max Planck Institute for Mathematics in the Sciences,
04103 Leipzig, Germany

\end{document}